\documentclass{amsart}

\usepackage{amsfonts}
\usepackage{amsmath}
\usepackage{amsthm}
\usepackage{amssymb}
\usepackage[english]{babel}
\usepackage{graphics}
\usepackage{latexsym}
\usepackage{longtable} 
\usepackage{mathrsfs} 
\usepackage{graphicx}
\usepackage{wrapfig} 
\usepackage[pdftex,colorlinks=true,linkcolor=blue,urlcolor=blue,citecolor=blue]{hyperref}

\usepackage{subfigure}

\newtheorem{theorem}{Theorem}
\newtheorem{corollary}[theorem]{Corollary}
\newtheorem{lemma}[theorem]{Lemma}
\newtheorem{proposition}[theorem]{Proposition}

\newtheorem{claim}[theorem]{Claim}
\newtheorem{example}[theorem]{Example}

\theoremstyle{definition}
\newtheorem{definition}[theorem]{Definition}
\newtheorem{remark}[theorem]{Remark}

\newcommand{\mA}{\mathcal{A}}
\newcommand{\mC}{\mathcal{C}}
\newcommand{\mL}{\mathcal{L}}
\newcommand{\mH}{\mathcal{H}}
\newcommand{\mF}{\mathcal{F}}

\newcommand{\mE}{\mathcal{E}}
\newcommand{\mN}{\mathcal{N}}
\newcommand{\mM}{\mathcal{M}}

\newcommand{\mK}{\mathcal{K}}

\newcommand{\mP}{\mathscr{P}}
\newcommand{\mO}{\mathcal{O}}
\newcommand{\mY}{\mathscr{Y}}
\newcommand{\mD}{\mathcal{D}}
\renewcommand{\S}{\mathcal{S}}

\newcommand{\A}{\mathrm{A}}
\newcommand{\R}{\mathbb{R}}

\newcommand{\N}{\mathbb{N}}

\newcommand{\mB}{\mathbb{B}}
\newcommand{\X}{\textbf{X}}

\renewcommand{\H}{\mathrm{H}}
\newcommand{\E}{\mathrm{E}}

\newcommand{\noi}{\noindent}
\newcommand{\ms}{\medskip}

\newcommand{\al}{\alpha}
\newcommand{\be}{\beta}
\newcommand{\ga}{\gamma}

\newcommand{\de}{\delta}
\newcommand{\De}{\Delta}
\newcommand{\e}{\varepsilon}
\newcommand{\si}{\sigma}

\newcommand{\la}{\lambda}

\newcommand{\Om}{\Omega}
\newcommand{\om}{\omega}

\newcommand{\ze}{\zeta}

\newcommand{\av}{-\!\!\!\!\!\!\int}

\newcommand{\D}{\mathrm{D}} 
 
\newcommand{\weak }{\, -\!\!\!\!-\!\!\!\rightharpoonup}
\newcommand{\weakstar }{ \overset{\, *_{\phantom{|}}}{{\smash{\weak }}\, } }

\newcommand{\larrow}{\longrightarrow}
\newcommand{\ot}{\otimes}

\newcommand{\ri}{\rightarrow}

\newcommand{\p}{\partial}
\newcommand{\sub}{\subseteq}
\newcommand{\set}{\setminus}
\newcommand{\by}{\times}

\newcommand{\sgn}{\mathrm{sgn}}
\newcommand{\ess}{\mathrm{ess}}

\newcommand{\dist}{\mathrm{dist}}

\newcommand{\supp}{\mathrm{supp}}

\newcommand{\bt}{\begin{theorem}}\newcommand{\et}{\end{theorem}}
\newcommand{\bd}{\begin{definition}}\newcommand{\ed}{\end{definition}}
\newcommand{\bl}{\begin{lemma}}\newcommand{\el}{\end{lemma}}
\newcommand{\beq}{\begin{equation}}\newcommand{\eeq}{\end{equation}}
\newcommand{\bc}{\begin{claim}}\newcommand{\ec}{\end{claim}}
\newcommand{\bex}{\begin{example}}\newcommand{\eex}{\end{example}}
\newcommand{\bcor}{\begin{corollary}}\newcommand{\ecor}{\end{corollary}}
\newcommand{\bp}{\begin{proof}}\newcommand{\ep}{\end{proof}}

\newcommand{\BPL}{\medskip \noindent \textbf{Proof of Lemma} }
\newcommand{\BPC}{\medskip \noindent \textbf{Proof of Claim} }
\newcommand{\BPCOR}{\medskip \noindent \textbf{Proof of Corollary} }
\newcommand{\BPP}{\medskip \noindent \textbf{Proof of Proposition} }
\newcommand{\BPT}{\medskip \noindent \textbf{Proof of Theorem} }

\numberwithin{equation}{section}

\begin{document}

\title[$2$nd Order $L^\infty$ Variational Problems and the $\infty$-Polylaplacian]{Second Order $L^\infty$ Variational Problems and the $\infty$-Polylaplacian}

\author{Nikos Katzourakis and Tristan Pryer}
  \thanks{\!\!\!\!\!\!\!\!\texttt{N.K. has been partially financially supported by the EPSRC grant EP/N017412/1. \\ T.P. has been partially financially supported by the EPSRC grant EP/P000835/11}}
    
\address{Department of Mathematics and Statistics, University of Reading, Whiteknights, PO Box 220, Reading RG6 6AX, UNITED KINGDOM}
\email{n.katzourakis@reading.ac.uk}

\email{t.pryer@reading.ac.uk}



\date{}



\begin{abstract} In this paper we initiate the study of $2$nd order variational problems in $L^\infty$, seeking to minimise the $L^\infty$ norm of a function of the hessian. We also derive and study the respective PDE arising as the analogue of the Euler-Lagrange equation. Given $\mathrm{H}\in C^1(\mathbb{R}^{n\times n}_s)$, for the functional
\[ \label{1}
\ \ \ \ \mathrm{E}_\infty(u,\mathcal{O})\, =\, \big\| \mathrm{H}\big(\mathrm{D}^2 u\big) \big\|_{L^\infty(\mathcal{O})}, \ \ \ u\in W^{2,\infty}(\Omega),\ \mathcal{O}\subseteq \Omega, \tag{1} 
\]
the associated equation is the fully nonlinear $3$rd order PDE
\[ \label{2}
\A^2_\infty  u\, :=\,\big(\mathrm{H}_X\big(\mathrm{D}^2u\big)\big)^{\otimes 3}:\big(\mathrm{D}^3u\big)^{\otimes 2}\, =\,0. \tag{2}
\]
Special cases arise when $\mathrm{H}$ is the Euclidean length of  either the full hessian or of the Laplacian, leading to the $\infty$-Polylaplacian and the $\infty$-Bilaplacian respectively. We establish several results for \eqref{1} and \eqref{2}, including existence of minimisers, of absolute minimisers and of ``critical point" generalised solutions, proving also variational characterisations and uniqueness. We also construct explicit generalised solutions and perform numerical experiments. 
\end{abstract}

\maketitle

\tableofcontents

\section{Introduction} \label{section1}

In this paper we initiate the study of higher order variational problems in the space $L^\infty$ and of their respective associated equations. As a first step, we consider the problem of minimising the $L^\infty$ norm of a function of the hessian and study its connection to a respective PDE arising as the analogue of the Euler-Lagrange equation. More precisely, let $\Om\sub \R^n$ be a bounded open set, $n\in \N$. For a real function $u\in C^3(\Om)$, the gradient, the hessian and the 3rd order derivative are denoted respectively by
\[
\begin{split}
\D u\, &=\, (\D_i u)_{i=1}^n\ \ \ \ \ \ \ : \ \Om\sub \R^n \larrow \R^n,
\\
\D^2 u \, &=\, \big( \D^2_{ij}u \big)_{i,j=1}^n\ \ \ \ : \ \Om\sub \R^n \larrow \R^{n^{\ot 2}}_s,
\\
\D^3 u \, &=\, \big( \D^3_{ijk}u \big)_{i,j,k=1}^n \ : \ \Om\sub \R^n \larrow \R^{n^{\ot 3}}_s,
\end{split}
\]
and in general for any $N\in \N$ the $N$th order derivative $\D^Nu$ is valued in the space $\R^{n^{\ot N}}_s$, the $N$-fold symmetric tensor power of $\R^n$:
\[
\R^{n^{\ot N}}_s\, :=\, \Big\{T \in \R^n \ot \cdots \ot \R^n \equiv \R^{n^{\ot N}}\ \Big| \ T_{i_1 ... i_N} =T_{\si(i_1 ... i_N)}, \, \si \text{ permutation}\Big\}.
\]
Given a fixed function $\H\in C^1\big( \R^{n^{\ot 2}}_s \big)$, we consider the supremal functional
\beq \label{1.1}
\begin{split}
\mathrm{E}_\infty(u,\mO)\, :=\, \left\| \mathrm{H}\big(\D^2 u\big) \right\|_{L^\infty(\mO)}, 
\ \ \ \ 
u\in W^{2,\infty}(\Om),\ \ \mO\sub \Omega \text{ measurable}.
\end{split}
\eeq
It turns out that the associated PDE which plays the role of ``Euler-Lagrange" equation for \eqref{1.1} is the following \emph{fully nonlinear PDE of $3$rd order}:
\beq \label{1.2}
\A^2_\infty  u\, :=\,\big(\mathrm{H}_X\big(\D^2u\big)\big)^{\ot 3}:\big(\D^3u\big)^{\ot 2}\, =\,0. 
\eeq
The notation $(\cdot)^{\ot N}$ symbolises the $N$-fold tensor power of the object in the bracket and $\H_X$ denotes the gradient of $\H$ with respect to its matrix argument, whilst ``$:$" is a higher order contraction which extends the usual Euclidean inner product of the space of matrices. In index form, \eqref{1.2} reads
\[
\sum_{i,j,k,l,p,q=1}^n \Big(\mathrm{H}_{X_{ij}}\big(\D^2u\big) \, \mathrm{H}_{X_{kl}}\big(\D^2u\big)\, \mathrm{H}_{X_{pq}}\big(\D^2u\big)\Big)\, \D^3_{ikl}u\, \D^3_{jpq}u\, =\, 0.
\]
Further, by contracting derivatives we may rewrite \eqref{1.2} as
\beq \label{1.2a}
\ \ \ \A^2_\infty  u\, =\, \mathrm{H}_X\big(\D^2u\big) : \, \D\big(\mathrm{H}\big(\D^2u\big)\big) \ot \D\big(\mathrm{H}\big(\D^2u\big)\big)\, =\,0.
\eeq
Special cases of the PDE \eqref{1.2} arise when $\H$ is either the Euclidean length of the hessian, $\D^2u$, in $\R^{n^{\ot 2}}_s$ (squared) $|\D^2u|^2=\sum_{i,j=1}^n \big(\D^2_{ij}u\big)^2 = \D^2u :\D^2u$ or the absolute value of the Laplacian, $\De u$ (squared), giving rise to the following respective PDEs which (borrowing the terminology from Differential Geometry, see e.g.\ \cite{M, HM, GS}) we call the \emph{$\infty$-Polylaplacian and the $\infty$-Bilaplacian}: 
\begin{align}  
\Pi^2_\infty u\, &:=\,\big( \D^2u \big)^{\ot 3}:\big(\D^3u\big)^{\ot 2}\, =\,0, \label{1.3}
\\
\De^2_\infty u\, & :=\,\big(\De u\, \mathrm{I}\big)^{\ot 3}:\big(\D^3u\big)^{\ot 2}\, =\,0 .
\label{1.4}
\end{align}
In \eqref{1.4}, ``$\mathrm{I}$" is the identity matrix. An equivalent way to write $\Pi^2_\infty$ and $\De^2_\infty$ (after a rescaling) is respectively
\[
\begin{split}
 \D^2u  : \, \D\big( |\D^2u |^2\big) \ot \D\big( |\D^2u |^2\big)\, =\,0,
\ \ \ \ 
(\De u)^3 \big|\D (\De u)\big|^2\, =\, 0.
\end{split}
\]
The study of $1$st order variational problems when minimising a function of the gradient 
\beq \label{1.5}
\ \ \ \ E_{\infty}(u,\mO)\, :=\, \left\|H(\D u) \right\|_{L^\infty(\mO)}, 
\ \ \ u\in W^{1,\infty}(\Om),\ \mO\sub \Omega \text{ measurable},
\eeq
is by now quite standard and has been pioneered by Aronsson who first considered \eqref{1.5} in the 1960s (\cite{A1}-\cite{A7}). In this case, the respective PDE is quasilinear of $2$nd order and is commonly known as the ``Aronsson equation":
\beq \label{1.6}
\A_\infty u\, :=\, H_p(\D u) \ot H_p(\D u) : \D^2u\, =\,0.
\eeq
An important special case of \eqref{1.6} arises for $H(p)=|p|^2$ and is known as the $\infty$-Laplacian: $\De_\infty u:=  \D u  \ot \D u : \D^2u =0$. The field has undergone a marvellous development since then, especially in the 1990s when the advent of the theory of Viscosity Solutions for fully nonlinear $2$nd order PDEs made possible the rigorous study of the non-divergence equation \eqref{1.6} and of its non-smooth solutions (for a pedagogical introduction with numerous references we refer to \cite{K7, C}). The popularity of this area owes to both the intrinsic mathematical interest as well as to the importance for applications, since minimisation of the maximum ``energy" provides more realistic models than the standard integral counterparts of average ``energy". Let us also note that the vectorial first order case is under active research and since the early 2010s is being developed very rapidly. To the best of our knowledge, the systematic study has been initiated by the first author (see \cite{K1}-\cite{K6}, \cite{K8}-\cite{K11}, and the joint contributions with Abugirda, Ayanbayev, Croce and Pisante \cite{AK,AyK,AAK,CKP}), while the second author has been working on the numerical analysis of them (\cite{KP, P, LP1, LP2}).

As it is well known from the $1$st order case of \eqref{1.5}, supremal functionals  lack ``locality" and the requirement of minimality has to be imposed at the outset on all subdomains, not just the domain itself as in the case of integral functionals. In particular, mere minimisers of \eqref{1.1} are not truly optimal and may not solve in any sense the PDE \eqref{1.2}. The variational principle we will be considering for \eqref{1.1} is the following extension of Aronsson's notion of Absolute Minimisers:

\begin{definition} \label{definition1} A function $u\in W^{2,\infty}(\Om)$ is called a \emph{$2$nd order Absolute Minimiser} of \eqref{1.1} when
\[
\ \ \ \ \E_\infty(u,\Om')\,\leq\, \E_\infty(u+\phi,\Om'),\ \ \ \forall\, \Om'\Subset \Om,\ \forall \phi\in W^{2,\infty}_0(\Om').
\]
The set of $2$nd order Absolute Minimisers will be symbolised by $\mathrm{AM}^2(\E_\infty,\Om)$.
\end{definition}
We would like to emphasise that, as the explorative results in this paper will make apparent,
\smallskip

\noi \emph{the higher order (scalar) case of \eqref{1.1}-\eqref{1.2} can not be developed ``by analogy" to the first order case (neither scalar nor vectorial) and unexpected phenomena arise}. 

\ms

For example, a fundamental difficulty associated to \eqref{1.1} is that \emph{even the $1$-dimensional problem of minimising $|u''|^2$ in $L^\infty$ (or in $L^p$) is not trivial; in particular, even in this case the minimisers are non-polynomial and actually have singular points, being non-$C^2$ and just $W^{2,\infty}$ (let alone $C^3$).} An extra difficulty associated to \eqref{1.4} is that the respective functional $(u,A) \mapsto \|(\De u)^2\|_{L^\infty(A)}$ is not coercive in $W^{2,\infty}(\Om)$ but instead in the space
\[
\bigg\{u \in \bigcap_{1<p<\infty}W^{2,p}(\Om)\ :\ \De u \in L^{\infty}(\Om)\bigg\}
\]
because of the failure of the Calderon-Zygmund $L^p$-estimates in the extreme case $p=\infty$ (see e.g.\ \cite{GM, GT}). More importantly and even more unexpectedly, the relevant PDE \eqref{1.2} is not any more $2$nd order quasilinear and degenerate elliptic, but instead \emph{$3$rd order fully nonlinear since it is quadratic in the highest order derivative; moreover, it is highly degenerate but in no obvious fashion elliptic}. To the best of our knowledge this is the first instance in Calculus of Variations in general where a fully nonlinear PDE of odd order is actually variational and in addition not a null Lagrangian. Even the $1$-dimensional version of the $\infty$-Polylaplacian/$\infty$-Bilaplacian is not trivial; in fact, for $n=1$ both equations \eqref{1.3}-\eqref{1.4} simplify to
\beq \label{1.9}
\De^2_\infty u\, =\, (u'')^3(u''')^2=\,0
\eeq
and it is easy to see that its \emph{solutions can not in general be $C^3$}. Accordingly, 
\[
\left\{
\begin{array}{rl}
\ \ (u'')^3(u''')^2=\,0, \ & \text{ in }\Om, \\
 \ \ u=g, \ u'=g', & \text{ on }\p\Om
 \end{array}
 \right.
\]
is \emph{solvable in the class of $C^3$ solutions if and only if $g$ is a quadratic polynomial}. In general, even the ``best candidate" solution coming from the limit of the $p$-Bilaplacian as $p\ri\infty$ \emph{is not $C^2$ and $u''$ is discontinuous}.

Perhaps the greatest difficulty associated to the study of \eqref{1.2} is that all standard approaches in order to define generalised solutions based on maximum principle or on integration-by-part considerations seem to fail. As highlighted above, there is a real necessity for such a notion for \eqref{1.2} even when $n=1$. More concretely, by a separation of variables of the form $u(x,y)=f(x)+g(y)$ on $\R^2$, one easily arrives as Aronsson did in \cite{A6} to the singular global $\infty$-Polyharmonic function on $\R^2$ $u(x,y)= |x|^{\frac{12}{5}}-\, |y|^{\frac{12}{5}}$ which is saddle-shaped but not thrice differentiable on the axes because $|\D^2u(x,y)|^2$ $\cong |x|^{4/5}+ |y|^{4/5}$. Further singular solutions without $3$rd order derivatives arise by the special class of solutions to the fully nonlinear $2$nd order equation $\H\big(\D^2u\big) =c$. 

Motivated in part by the systems arising in vectorial Calculus of Variations in $L^\infty$, the first author has recently introduced in \cite{K8, K9} a new efficient theory of generalised solutions which applies to fully nonlinear systems of any order
\[ 
\ \ \  \mF\Big(\cdot,u ,\D u ,\D^2u,...,\D^N u\Big)\, =\, 0, \quad \text{ in }\Om,
\]
and allows for \emph{merely measurable} mappings as solutions. This general approach of the so-called \emph{$\mD$-solutions} is based on the probabilistic representation of those derivatives which do not exist classically. The tool in achieving this is the weak* compactness of difference quotients in the Young measures valued into a compactification of the ``space of jets". For the special case of the $3$nd order PDE \eqref{1.2}, we can motivate the idea as follows: let $u$ be a $W^{3,\infty}(\Om)$ strong solution of
\beq   \label{1.10}
\ \ \ \mF\big(\D^2 u,\D^3u\big)\,=\, 0, \quad \text{a.e.\ on }\Om.
\eeq
In order to interpret the $3$rd derivative rigorously for just $W^{2,\infty}_{\text{loc}}(\Om,\R^N)$ (which is the natural regularity class for \eqref{1.2} arising from \eqref{1.1}), we argue as follows: let us restate \eqref{1.10} as
\beq  \label{1.11}
\ \ \ \int_{\R^{n^{\ot 3}}_s} \Phi(\X)\, \mF\big(\D^2 u(x),\X\big)\, d[\de_{\D^3 u(x)} ](\X)\, =\, 0, \quad \text{ a.e. }x\in \Om,
\eeq
for any $\Phi \in C_c\big( \R^{n^{\ot 3}}_s \big)$ with compact support. Namely, we view the $3$rd derivative tensor $\D^3 u$ as a probability-valued mapping $\Om\sub\R^n \larrow \mathscr{P}\big(\R^{n^{\ot 3}}_s \big)$ which is given by $x\mapsto \de_{\D^3 u(x)}$, the Dirac measure at the $3$rd derivative. Also, we may rephrase that $\D^3 u$ is the a.e.\ sequential limit of the difference quotients $\D^{1,h}\D^2 u$ of the hessian along infinitesimal sequences $(h_m)_1^\infty$ by writing
\beq  \label{1.12}
\ \ \ \de_{\D^{1,h_m}\D^2 u} \weakstar\, \de_{\D^3 u}, \quad \text{as }m\ri\infty.
\eeq
The weak* convergence in \eqref{1.12} is taken in the set of Young measures valued into the tensor space $\R^{n^{\ot 3}}_s $ (the set of weakly* measurable probability-valued maps $\Om\sub\R^n \larrow \mathscr{P}\big(\R^{n^{\ot 3}}_s \big)$, for details see Section \ref{section2} and \cite{CFV, FG, V, Pe, FL}). The idea arising from \eqref{1.11}-\eqref{1.12} is that perhaps general probability-valued ``diffuse $3$rd derivatives" could arise for twice differentiable maps which may not be the concentrations $\de_{\D^3 u}$. This is actually possible upon replacing $\R^{n^{\ot 3}}_s $ by its $1$-point sphere compactification in order to gain some compactness: $\smash{\overline{\R}}^{n^{\ot 3}}_s := \R^{n^{\ot 3}}_s\cup\{\infty\}$. Then, the maps $(\de_{\D^{1,h}\D^2 u})_{h\neq 0}$ considered as Young measures do have subsequential weak* limits which can play the role of generalised $3$rd order derivatives:

\begin{definition}[Diffuse $3$rd derivatives, cf.\ \cite{K8,K9}]  \label{definition2}
For any $ u \in W^{2,1}_{\text{loc}}(\Om)$, we define its diffuse $3$rd derivatives $\mD^3 u$ as the limits of difference quotients of $\D ^2u$ along infinitesimal sequences $(h_m)_1^\infty$ in the Young measures valued into $\smash{\overline{\R}}^{n^{\ot 3}}_s$:
\[
\ \ \  \ \de_{\D^{1,h_m}\D^2 u}\weakstar \mD^3 u, \  \ \ \, \text{ in }\mY\big(\Om,\smash{\overline{\R}}^{n^{\ot 3}}_s\big), \ \ \ \text{ as }m\ri \infty.
\]
If $\{e^1,...,e^n\}$ stands for the standard basis of $\R^n$, then apparently
\[
\D^{1,h}v\, :=\,  \big( \D^{1,h}_1,..., \D^{1,h}_n v\big),\ \ \D^{1,h}_iv(x)\, :=\, \frac{1}{h}\big[v(x+he^i)-v(x)\big] ,\ \ h\neq 0.
\] 
\end{definition} 
Since the set $\mY\big(\Om,\smash{\overline{\R}}^{n^{\ot 3}}_s\big)$ is weakly* compact, every map possesses at least one diffuse $3$rd derivative and actually exactly one if the hessian is a.e.\ differentiable with measurable derivative (see \cite{K8}).

\begin{definition}[Twice differentiable $\mD$-solutions of $3$rd order PDEs, cf.\ \cite{K8}] \label{definition3} Let $\mF: \R^{n^{\ot 2}}_s\!\by \R^{n^{\ot 3}}_s \larrow \R$ be Borel measurable. A function  $u \in W^{2,1}_{\text{loc}}(\Om,\R^N)$ is a $\mD$-solution to 
\beq \label{2.11}
\mF\big(\D^2 u,\D^3 u\big)\, =\,0, \ \ \text{ in }\Om,
\eeq
when for any diffuse $3$rd derivative $\mD^3 u \in \mY\big(\Om,\smash{\overline{\R}}^{n^{\ot 3}}_s\big)$ and any $\Phi \in C_c\big( {\R}^{n^{\ot 3}}_s \big)$,
\[
\ \ \ \int_{\smash{\overline{\R}}^{n^{\ot 3}}_s} \Phi(\X)\, \mF\big(\D^2 u(x),\X \big)\, d\big[\mD^3 u(x) \big](\X)\, =\, 0, \ \ \text{ a.e.\ }x\in \Om.
\] 
 \end{definition} 
The notion of generalised solution of Definitions \ref{definition2} \& \ref{definition3} will be the central notion of solution for our fully nonlinear PDE \eqref{1.2}. For more on the theory of $\mD$-solutions for general systems, analytic properties, existence/uniqueness/partial regularity results see \cite{K8}-\cite{K11} and \cite{AyK,AAK,AK,CKP}.

We note that the interpretation of \eqref{1.2} in the ``contracted" form \eqref{1.2a} is not generally appropriate for non-$C^3$ solutions; interpreting the ``expanded" equation \eqref{1.2} in a weak sense is essential. In particular, even in the $1D$ case of \eqref{1.9}, the results of Sections \ref{section8}-\ref{section9} demonstrate that seeing the $\infty$-Poly/Bilaplacian as $ u''\big(\big(|u''|^2\big)' \big)^2$ is \emph{not appropriate even when $n=1$ since there exist solutions for which $u''$ is piecewise constant and  hence the distributional derivative $\big(|u''|^2\big)'$ is a measure} (whose square is not well defined!).

In this paper we are concerned with the study of $2$nd order Absolute Minimisers of \eqref{1.1}, of $\mD$-solutions to \eqref{1.2} and of their analytic properties and their connection. To this end, we prove several things and the table of contents is relatively self-explanatory of the results we obtain in this paper. Below we give a quick description and main highlights:

In Section \ref{section2} we give a quick review of the very few ingredients of Young measures into spheres which are utilised in this paper for the convenience of the reader.

In Section \ref{section3} we formally derive the equation \eqref{1.1} in the limit of the Euler-Lagrange equation of the respective $L^p$ functional $u \mapsto \|\H\big(\D^2u\big)\|_{L^p(\Om)}$ as $p\ri \infty$.

In Section \ref{section5} we characterise $C^3$ solutions to \eqref{1.2} via the flow map of an ODE system along the orbits of which the energy is constant (Proposition \ref{proposition4}).

In Section \ref{section4} we prove existence of minimisers for \eqref{1.1} given Dirichlet boundary condition on a bounded open set under two sets of weak hypotheses which include both the $\infty$-Polylaplacian \eqref{1.3} and the $\infty$-Bilaplacian \eqref{1.4} (Theorems \ref{theorem4}, \ref{theorem4a}). We also give a complete solution to the problem of existence-uniqueness and description of the fine structure of $2$nd order Absolute Minimisers of \eqref{1.1} and of the corresponding minimising $\mD$-solutions to \eqref{1.2} when $n=1$ (Theorem \ref{theorem6}).

In Section \ref{section6} we establish the necessity of the PDE \eqref{1.2} for $2$nd order Absolute Minimisers of \eqref{1.1} in the class of $C^3$ solutions (Theorem \ref{theorem5}(I)). This is nontrivial even for $C^3$ solutions because standard $1$st order arguments fail to construct test functions in $W^{2,\infty}_0(\Om)$ and a deep tool is required, the Whitney extension theorem (\cite{W, M, F}). If further $\H$ depends on $\D^2u$ via the projection $A\!:\!\D^2u$ along a fixed matrix (e.g.\ on the Laplacian $\De u=\D^2u \!:\!\mathrm{I}$), we prove sufficiency as well and hence equivalence (Theorem \ref{theorem5}(II)). As a consequence, in the latter case we deduce uniqueness in the $C^3$ class for \eqref{1.1} and \eqref{1.2} (Corollary \ref{corollary6}).

In Section \ref{section7} we employ the Dacorogna-Marcellini Baire Category method (\cite{DM, D}) which in a sense is the analytic counterpart to Gromov's Convex integration and establish the existence of non-minimising ``critical point" $\mD$-solutions to the Dirichlet problem \eqref{1.2} (Theorem \ref{theorem12a}). We construct $\mD$-solutions in $W^{2,\infty}_g(\Om)$ with the extra geometric property of solving strongly the fully nonlinear equation $\H\big(\D^2u\big)=c$ at any large enough energy level $c>0$. Interestingly, \emph{we do not assume any kind of convexity or level-convexity or ``BJW-convexity"} (i.e.\ the notion of $L^\infty$-quasiconvexity introduced in \cite{BJW2}). This method has previously been applied in the construction of critical point vectorial $\mD$-solutions in the first order case in \cite{K9, CKP} and has some vague relevance to the method used in \cite{DS} to construct solutions to the Euler equations.

In Section \ref{section8} we solve explicitly the $p$-Bilaplacian (weakly) and the $\infty$-Bilaplacian (in the $\mD$-sense) in the case $n=1$ (Theorem \ref{theorem12}). In particular, the $p$-Biharmonic functions are $C^\infty$ except for at most one point in the domain and $\infty$-Biharmonic functions are smooth except for at most two points in the domain. In the $\infty$-case, the $\mD$-solutions we construct are non-minimising but have fixed energy level and this allows to have uniqueness (absolutely minimising $\mD$-solutions are constructed in Theorem \ref{theorem6}).

Finally in Section \ref{section9} we perform some numerical experiments by considering the solutions of the $p$-Bilaplacian for given Dirichlet data and large $p$. The experiments confirm numerically that the $2$nd derivatives of $\infty$-Biharmonic functions generally can not be continuous (this is actually proved in Theorems \ref{theorem6} and \ref{theorem5}). Even more interestingly, for $n=2$ the Laplacian of $\infty$-Biharmonic functions appears to be piecewise constant and we have the emergence of non-trivial interfaces whereon it is discontinuous and actually changes sign when crossing the interface. Further, for the ``balanced" symmetric energy $\H(X)=|X \!:\!\mathrm{I}|^2$ we are using, the values of Laplacian appear to be opposite on the phases, whilst the absolute value of it seems to extend to a constant function throughout the domain (this is proved in Theorem \ref{theorem6} for $n=1$).

Although in this paper we do not consider any immediate applications of our results, we would like to point out that $2$nd order minimisation problems in $L^\infty$ are very important in several areas of pure and applied Mathematics. In particular, in the papers \cite{MS} and \cite{S} the authors consider the problem of minimising in $L^\infty$ the Gaussian curvature (if $n=2$) and the scalar curvature (if $n\geq3$) of a fixed background Riemannian manifold over a conformal class of deformations of the metric. Although they consider (mere) minimisers and not absolute minimisers of their geometric functionals which is the appropriate notion in $L^\infty$, the method they use to construct them is via $L^p$ approximations and this seems to select the ``good" absolutely minimising object. The theoretical and numerical observations we make in this paper are compatible with phenomena of piecewise constant energy and interfaces of discontinuities in the differential-geometric context of minimisation of the curvature in $L^\infty$ in \cite{MS,S}. Also, after this work had been completed and appeared as a preprint, we learned that in the paper \cite{AB} Aronsson and Barron had already previously derived our PDE \eqref{2} (Remark 4.9 on p.\ 78) without investigating it further. We conclude by noting that in our companion paper \cite{KP2} we establish rigorous numerical approximations for $\infty$-Polyharmonic and $\infty$-Biharmonic functions and we also consider concrete applications.

\section{A quick guide of Young measures valued into spheres} \label{section2}

Here we collect some rudiments of Young measures taken from \cite{K8} which can be found in greater generality and different guises e.g.\ in \cite{CFV, FG}. Let $\Om\sub \R^n$ be open and let us consider the Banach space of $L^1$ maps valued in the space of continuous functions over the compact manifold $\smash{\overline{\R}}^{n^{\ot 3}}_s$ (for more see e.g.\ \cite{Ed, FL}): $L^1\big( \Om, C\big(\smash{\overline{\R}}^{n^{\ot 3}}_s\big)\big)$. The elements of this space are Carath\'eodory functions $\Phi : \Om \by \smash{\overline{\R}}^{n^{\ot 3}}_s \larrow \R$ satisfying 
\[
\| \Phi \|_{L^1( \Om, C(\R^{n^{\ot 3}}_s))}\, =\, \int_\Om \bigg\{\max_{\X \in \smash{\overline{\R}}^{n^{\ot 3}}_s}\big|\Phi(x,\X)\big|\bigg\}\, dx\, <\, \infty.
\]
The dual of this Banach space is given by
\[
L^\infty_{w^*}\big( \Om, \mM \big(\smash{\overline{\R}}^{n^{\ot 3}}_s\big)\big)\, =\, L^1\big( \Om, C\big(\smash{\overline{\R}}^{n^{\ot 3}}_s\big)\big)^* 
\]  
and consists of measure-valued maps $x \mapsto \vartheta(x) $ which are weakly* measurable, namely the real function $x\mapsto [\vartheta(x)](\mathcal{U})$  is measurable on $\Om$ for any fixed open set $\mathcal{U} \sub \smash{\overline{\R}}^{n^{\ot 3}}_s$. The unit closed ball of this dual space is sequentially weakly* compact because the $L^1$ space above is separable. The duality pairing is given by
\beq \label{dp}
\left\{
\begin{split}
\ \ \ \ \langle\cdot,\cdot\rangle\ :\ \ & \ L^\infty_{w^*}\big( \Om, \mM\big(\smash{\overline{\R}}^{n^{\ot 3}}_s\big)\big) \by L^1\big( \Om, C\big(\smash{\overline{\R}}^{n^{\ot 3}}_s\big)\big) \larrow \R,  \ \ 
\\
& \langle \vartheta, \Phi \rangle\, :=\, \int_\Om \int_{\smash{\overline{\R}}^{n^{\ot 3}}_s} \Phi(x,\X)\, d[\vartheta(x)] (\X)\, dx.
\end{split}
\right.
\eeq

\noi \textbf{Definition} (Young Measures). {\it The subset of the unit sphere consisting of probability-valued mappings comprises the Young measures from $\Om\sub \R^n$ to the sphere $\smash{\overline{\R}}^{n^{\ot 3}}_s$:}
\[
\mY\big(\Om,\smash{\overline{\R}}^{n^{\ot 3}}_s\big)\, :=\, \Big\{ \vartheta \in L^\infty_{w^*}\big( \Om, \mM\big(\smash{\overline{\R}}^{n^{\ot 3}}_s\big)\big)\, : \, \text{ a.e. }x\in \Om,\ \vartheta(x) \in \mP\big(\smash{\overline{\R}}^{n^{\ot 3}}_s\big)\Big\}.
\]
We now record for later use the next standard facts (for the proofs see e.g.\ \cite{FG}): 
\smallskip

\noi (i) \emph{Any measurable map $U : \Om\sub \R^n \larrow \R^{n^{\ot 3}}_s$ induces a Young measure $\de_U$ given by $\de_U(x):= \de_{U(x)}$, $x\in \Om$.} 

\noi (ii) \emph{The set $\mY\big(\Om,\smash{\overline{\R}}^{n^{\ot 3}}_s\big)$ is sequentially weakly* compact and convex. In particular, any sequence $(\de_{U^m})_1^\infty$ has a subsequence such that $\de_{U^{m_k}} \weakstar \vartheta$ as $k\ri\infty$.}  

\smallskip

\noi (iii) \emph{If $(U^m)_1^\infty,U^\infty$ are measurable maps $\Om\sub \R^n\larrow \R^{n^{\ot 3}}_s$, then $\de_{U^{m}} \weakstar \de_{U^\infty}$ in $\mY\big(\Om, \smash{\overline{\R}}^{n^{\ot 3}}_s\big)$ iff $U^m \larrow U^\infty$ a.e. on $\Om$, after perhaps the passage to subsequences.}

\smallskip

\noi (iv) \emph{The following is an one-sided characterisation of weak* convergence: 
$\vartheta^m\weakstar \vartheta$ as $m\ri \infty$ in $\mY\big(\Om, \smash{\overline{\R}}^{n^{\ot 3}}_s\big)$ iff $\langle \vartheta,\Psi\rangle \leq {\lim \inf}_{m\ri \infty}\langle \vartheta^m,\Psi\rangle$ for any function $\Psi:\Om \by \smash{\overline{\R}}^{n^{\ot 3}}_s \larrow (-\infty,+\infty]$ bounded from below, measurable in $x$ for all $\X$ and lower semicontinuous in $\X$ for a.e.\ $x$.}

\section{Derivation of the $L^\infty$ equation from $L^p$ as $p\ri\infty$}  \label{section3}
 
In this section we \emph{formally derive} the equation \eqref{1.2} in the limit of the Euler-Lagrange equations of the respective $L^p$ functionals
\beq \label{3.1}
\E_p(u,\Om)\,:=\, \left(\av_\Om \H\big(\D^2u(x)\big)^p\,dx\right)^{1/p}
\eeq
as $p\ri \infty$. Here the bar denotes average. The idea of approximating an $L^\infty$ variational problem by $L^p$ problems is quite standard by now for $1$st order problems in both the scalar and the vectorial case (see e.g.\ \cite{C, K7, K9, P, KP}) and has borne substantial fruit. Heuristically, this expectation stems from the fact that for a fixed function $u\in W^{2,\infty}(\Om)$, we have $\E_p(u,\Om) \larrow \E_\infty(u,\Om)$
 as $p\ri \infty$. For \eqref{3.1}, the Euler-Lagrange equation is the $4$th order divergence structure PDE
\beq \label{3.2}
\D^2:\Big(\H\big(\D^2u\big)^{p-1}\H_X\big(\D^2u\big)\Big)\,=\,0
\eeq
which in index form reads $\sum_{i,j=1}^n \D_{ij}^2\big(\H\big(\D^2u\big)^{p-1}\H_{X_{ij}}\big(\D^2u\big)\big)=0$. By distributing derivatives and rescaling, a calculation gives
\[
\begin{split}
 \sum_{i,j,k,l,p,q}  \Big(   \mathrm{H}_{X_{ij}}  & \big(\D^2u\big) \, \mathrm{H}_{X_{kl}}\big(\D^2u\big)\, \mathrm{H}_{X_{pq}}\big(\D^2u\big)\Big)\, \D^3_{ikl}u\, \D^3_{jpq}u
 \\
 =\, &-\, \frac{\H\big(\D^2u\big)}{p-2} \Bigg\{
 \sum_{i,j,k,l} \bigg(
\D_i\big( \H_{X_{kl}}\big(\D^2u\big)\big)\, \D^3_{jkl} \, \H_{X_{ij}}\big(\D^2u\big) \hspace{40pt}
\end{split}
\]
\[
\begin{split}
&\, +\, \D_i\big( \H\big(\D^2u\big)\big)\,  \D_j\big(\H_{X_{ij}}\big(\D^2u\big)\big)
\, +\, \D_j\big( \H\big(\D^2u\big)\big)\, \D_j\big(\H_{X_{ij}}\big(\D^2u\big)\big)
\\
&\, +\, \H_{X_{ij}}\big(\D^2u\big)\, \H_{X_{kl}}\big(\D^2u\big)\, \D^4_{ijkl}u
\bigg)\Bigg\}
\\
\ \ \ \ \ \ & -\, \frac{\H\big(\D^2u\big)^2}{(p-1)(p-2)} 
\Bigg\{
 \sum_{i,j,k,l,p,q}\bigg( 
 \H_{X_{ij}X_{kl}X_{pq}}\big(\D^2u\big) \, \D^3_{ipq}u \, \D^3_{jkl}u\, 
 \\
 &+\,  \H_{X_{ij}X_{kl}}\big(\D^2u\big) \, \D^4_{ijkl}u\bigg)\Bigg\}.
\end{split}
\]
Hence, we obtain $\A^2_\infty u = \big(\H_X\big(\D^2u\big)\big)^{\ot 3} : \big( \D^3 u\big)^{\ot 2} = o(1)$ as $p\ri \infty$ and hence we obtain \eqref{1.2}. Note that the Euler-Lagrange equation of \eqref{3.1} is $4$th order and quasilinear, while the limiting equation is a highly degenerate $3$rd order \emph{fully nonlinear} equation. In the next sections we utilise this device of $L^p$ approximations and prove rigorously the existence of minimisers and absolute minimisers.

\section{Characterisation of $\A^2_\infty $ via the flow map of an ODE}  \label{section5}
 
In this brief section, inspired by the $1$st order case (see \cite{C, K7, K1}) we give a description of classical solutions to our fully nonlinear PDE \eqref{1.1} in terms of the flow of a certain ODE system. In the $1$st order case the relevant ODE is a gradient system but in the present case it is more complicated and involves $3$rd order derivatives.
 
\begin{proposition} \label{proposition4} Let $\H\in C^1\big( \R^{n^{\ot 2}}_s \big)$, $\Om\sub \R^n$ an open set and $u\in C^3(\Om)$. Consider the continuous vector field
\beq \label{5.1}
\mathscr{V}\,:=\, \H_X\big(\D^2u\big) \D\big(\H\big(\D^2u\big)\big)\ :\ \ \Om\sub \R^n \larrow \R^n
\eeq
and the initial value problem
\beq \label{5.2}
\left\{
\begin{array}{ll}
\dot{\ga}(t)\, =\, \sgn\big(\mathscr{V}\big(\ga(t)\big)\big), &\ \ t\neq 0,\ms
\\
 \ga(0)\, =\,x, &\ \ t=0,
\end{array}
\right.
\eeq
where the initial condition is noncritical, i.e.\ $\mathscr{V}(x)\neq0$ and ``sgn" symbolises the sign. Then, along the trajectory we have the differential identity
\beq \label{5.3}
\big|\mathscr{V}\big(\ga(t)\big)\big|\frac{d}{dt} \Big(\H\big(\D^2u\big(\ga(t)\big)\big)\Big)\, =\, \A^2_\infty  u\big(\ga(t)\big) 
\eeq
and hence 
\[
\text{$\A^2_\infty  u=0$ in $\Om$} \ \ \ \Longleftrightarrow \ \ \  
\left\{
\begin{array}{l}
\forall\, x\in \Om\, :\, \mathscr{V}(x)\neq0, \ \text{there is a $C^1$} \smallskip
\\
\text{solution $\ga : (-\e,\e) \larrow \Om$ of \eqref{5.2}:}  \smallskip
\\
\H\big(\D^2u\big(\ga(t)\big)\big)\! \equiv \H\big(\D^2u(x)\big), \, |t|<\e.
\end{array}
\right.
\]
\end{proposition}
Note that, unlike the counterpart $1$st order case, the solution of initial value problem \eqref{5.2} may not be unique in general. 

\BPP \ref{proposition4}. In order to conclude it suffices to establish \eqref{5.3} the proof of which is a straightforward calculation. The proposition ensues.       \qed
\ms

Heuristically, the meaning of this result is the following: in view of \eqref{1.2} the functions with $\H\big(\D^2u\big)\equiv c$ are special solutions to the PDE. Conversely, all solutions satisfy $\H\big(\D^2u\big)\equiv c$ at least locally along the trajectories of the $1$st order ODE system \eqref{5.1}-\eqref{5.2}.

\section{Existence of $2$nd order Minimisers \& Absolute Minimisers}  \label{section4}

Herein we consider the problem of existence of minimisers and of $2$nd order Absolute Minimisers for \eqref{1.1} with given boundary values (Definition \ref{definition1}). To this end we will assume that $\H$ is \emph{level-convex} (namely has convex sub-level sets) and we will obtain our $L^\infty$ objects in the limit of \emph{approximate} minimisers of $L^p$ functionals. The methods of this section have been inspired by the paper \cite{BJW1} wherein the authors prove the existence of absolute minimisers when the rank of the gradient is at most one (scalar-valued functions or curves, see also the papers \cite{K8, AK} for relevant ideas). We begin below with the simpler case of the existence of (mere) minimisers and subsequently we will show that the candidate we construct is indeed in $\mathrm{AM}^2(\E_\infty,\Om)$ when $n=1$. 

\begin{theorem}[Existence of $L^\infty$ Minimisers and their $L^p$-approximation, I] \label{theorem4}
Let $n\in \N$, $\Om\sub \R^n$ bounded and open and $\H \in C(\R^{n^{\ot 2}}_s) $ a non-negative level-convex function (that is for any $t\geq 0$, the set of matrices $\{\H\leq t\}$ is convex in $\smash{\R^{n^{\ot 2}}_s}$). Suppose also there exist $C_1,C_2,r>0$ such that
\[
\ \ \ \ \H(X)\,\geq\, C_1|X|^r -\, C_2, \ \  \ \ X\in \R^{n^{\ot 2}}_s.
\]
Then, for any $g\in W^{2,\infty}(\Om)$ there exists a function $u_\infty\in W_g^{2,\infty}(\Om)$ such that:

\smallskip

\noi \emph{(a)} $u_\infty$ is a minimiser of \eqref{1.1} on $\Om$, i.e.\ $\E_\infty(u,\Om) \leq \E_\infty(\psi,\Om)$, for any $\psi\in W_g^{2,\infty}(\Om)$. 
\smallskip

\noi \emph{(b)} For any $q\geq 1$, $u_\infty$ is the weak $W^{2,q}(\Om)$-limit of a sequence of approximate minimisers $(u_p)_{p=1}^\infty$ of the integral functionals \eqref{3.1} placed in $W_g^{2,\infty}(\Om)$ along a subsequence $p_j\ri \infty$. Namely, for any $q\geq 1$ we have $u_{p_j}\weak u_\infty$ in $W^{2,q}(\Om)$ as $j\ri \infty$ and $u_p$ satisfies
\[
\ \ \ \ \E_p(u_p,\Om)^p\, \leq\ 2^{-p^2} +\, \inf\Big\{\E_p(\cdot,\Om)^p\, :\ W_g^{2,\infty}(\Om)\Big\}.
\]
\noi \emph{(c)} For any measurable $A\sub \Om$, we have the ``diagonal lower semi-continuity"
\[
\E_\infty(u_\infty,A)\, \leq\, \liminf_{j\ri \infty}\, \E_{p_j}(u_{p_j},A).
\]
\end{theorem}
The proof of Theorem \ref{theorem4} can be done mutatis mutandis to the proof of 
Theorem \ref{theorem4a} that follows and hence we refrain from giving the details. Theorem \ref{theorem4} does not include the case of the $\infty$-Bilaplacian \eqref{1.4} when minimising $(\De u)^2$ and more generally when $\H(X)=\textbf{H}(A:X)$ for some matrix $A>0$. In this case the appropriate space to obtain existence of minimisers is not $W^{2,\infty}(\Om)$ but instead the larger space 
\beq \label{space}
\mathcal{W}^{2,\infty}(\Om)\,:=\, \bigg\{u \in \bigcap_{1<p<\infty}W^{2,p}(\Om)\ \Big|\  \ A:\D^2u \in L^{\infty}(\Om)\bigg\}
\eeq
(with $\mathcal{W}_0^{2,\infty}(\Om)$ being defined in the obvious way) because of the inability to estimate $\D^2u$ in terms of $A\!:\!\D^2u$ in the $L^\infty$ norm (see \cite{GM}).

\begin{theorem}[Existence of $L^\infty$ Minimisers and their $L^p$-approximation, II] \label{theorem4a}
Let $n\in \N$, $\Om\sub \R^n$ a bounded open set and $\textbf{\emph{H}} \in C(\R)$ a non-negative level-convex function (i.e.\ for any $t\geq 0$, the sets $\{\textbf{\emph{H}}\leq t\}$ are intervals). Suppose also there exist $C_1,C_2,r>0$ such that
\[
\ \ \ \ \textbf{\emph{H}}(t)\, \geq\, C_1|t|^r -\, C_2, \ \ \ \ t\in \R.
\]
Let also $A\in \R^{n^{\ot 2}}_s$ be a strictly positive matrix.  Then, for any $g\in \mathcal{W}^{2,\infty}(\Om)$ (see \eqref{space}) there exists a function $u_\infty\in \mathcal{W}_g^{2,\infty}(\Om)$ such that:

\smallskip

\noi \emph{(a)} $u_\infty$ is a minimiser of the functional \eqref{3.1} for $\H(X)=\textbf{\emph{H}}(A\!:\!X)$ over the space $\mathcal{W}_g^{2,\infty}(\Om)$, i.e.\ $\E_\infty(u,\Om) \leq \E_\infty(\psi,\Om)$ for any $\psi\in \mathcal{W}_g^{2,\infty}(\Om)$.
\smallskip

\noi \emph{(b)} For any $q\geq 1$, $u_\infty$ is the weak $W^{2,q}(\Om)$-limit of a sequence of approximate minimisers $(u_p)_{p=1}^\infty$ of the integral functionals \eqref{3.1} (with $\H(X)=\textbf{\emph{H}}(A\!:\!X)$) placed in $\mathcal{W}_g^{2,\infty}(\Om)$ along a subsequence $p_j\ri \infty$. That is, for any $q\geq 1$  $u_{p_j}\weak u_\infty$ in $W^{2,q}(\Om)$ as $j\ri \infty$, whilst $u_p$ satisfies
\[
\ \ \ \ \E_p(u_p,\Om)^p\, \leq\ 2^{-p^2} +\, \inf\Big\{\E_p(\cdot,\Om)^p\, :\ \mathcal{W}_g^{2,\infty}(\Om)\Big\}.
\]
\noi \emph{(c)} The same lower semi-continuity statement as in Theorem \ref{theorem4}(c) holds true  but for the functionals \eqref{1.1} and \eqref{3.1} with $\H(X)=\textbf{\emph{H}}(A\!:\!X)$.
\end{theorem}

The idea of the proof of Theorem \ref{theorem4a} follows similar lines to those of the $1$st order results of \cite[Theorem 2.1]{BJW1}, \cite[Lemma 4]{K8}, \cite[Lemma 5.1]{K9}. As in \cite{BJW1}, the essential point is the use of Young measures (valued in the Euclidean space, in contrast to the sphere-valued Young measure we employ in the definition of $\mD$-solutions) in order to circumvent the lack of quasi-convexity for the $L^p$ approximating functionals for which the infimum may not be attained at a minimiser (hence the need for approximate minimisers at the $L^p$ level).

\BPT \ref{theorem4a}. We begin by noting that our coercivity lower bound and H\"older inequality imply the estimate
\beq \label{4.2}
\frac{C_1}{2}\left(\av_\Om \big|A:\D^2v\big|^{k}\right)^{\frac{r}{k}}-\, C_2 \, \leq\, \left(\av_\Om \textbf{H}\big(A:\D^2v\big)^{p}\right)^{\frac{1}{p}} \,\leq\, \E_\infty(v,\Om)
\eeq
for any $v\in \mathcal{W}^{2,\infty}_g(\Om)$ and $k\leq rp$. Fix $p>1+(1/r)$ and consider a minimising sequence $(u_{p,i})_{i=1}^\infty \sub \mathcal{W}^{2,\infty}_g(\Om)$ of $\E_p(\cdot,\Om)^p$. Select $i=i(p)$ large so that
\beq  \label{4.2a}
 \E_p(u_p,\Om)^p \, \leq\, 2^{-p^2}\, +\, \inf\Big\{\E_p(\cdot,\Om)^p\, :\ \mathcal{W}_g^{2,\infty}(\Om)\Big\}
\eeq
where $u_p:=u_{p,i(p)}$. In particular, $\E_p(u_p,\Om)\leq \E_\infty(g,\Om)+1$ and by \eqref{4.2} the functions $(A\!:\!\D^2u_p)_1^\infty$ are bounded in $L^k(\Om)$ for any fixed $k\in \N$. Since $A$ is constant and $u_p-g \in W_0^{2,k}(\Om)$, by the Calderon-Zygmund $L^k$-estimates (se e.g.\ \cite{GT}) and Poincar\'e inequality, we have
\[
\big\| A: \D^2u_p- A:\D^2g\big\|_{L^k(\Om)} \,\geq\ C(k,n,A)\, \| u_p-g \|_{W^{2,k}(\Om)}.
\]
Thus, the sequence $(u_p)_{p=1}^\infty$ is weakly precompact in $W^{2,k}(\Om)$ for any fixed $k\in \N$ and there exists $u_\infty \in \bigcap_{1<k<\infty} W^{2,k}(\Om)$ such that $u_p \weak u_\infty$ as $p\ri\infty$ along perhaps a subsequence. Further, by setting $v=u_p$ in \eqref{4.2}, using the weak lower-semicontinuity of the $L^k(\Om)$-norm and letting $p\ri \infty$ and $k\ri\infty$, we obtain $A\!:\!\D^2u_\infty \in L^\infty (\Om)$ which allows to infer that $u_\infty \in \mathcal{W}_g^{2,\infty}(\Om)$. Further, by \eqref{4.2a} and H\"older inequality, for any $\psi \in \mathcal{W}_g^{2,\infty}(\Om)$ and any $q\leq p$ we have
 \beq \label{4.3}
\left(\av_\Om \textbf{H}\big(A:\D^2u_p\big)^{q}\right)^{\frac{1}{q}} \leq\ 2^{-p}\,+\, \E_\infty(\psi,\Om) .
\eeq 
Consider now the sequence of Young measures generated by the scalar functions $(A\!:\!\D^2u_p)_1^\infty$, that is  we consider $(\de_{A:\D^2u_p})_{p=1}^\infty \sub \mY\big(\Om,\overline{\R}\big)$. Then, along perhaps a further subsequence we have $\de_{A:\D^2u_p} \weakstar \vartheta_\infty$ as $p\ri\infty$. Since $A\!:\!\D^2u_\infty \in L^\infty (\Om)$, the Young measure $\vartheta_\infty$ has compact support in $\R$. Thus, there is $R>0$ such that $\supp\big(\vartheta_\infty(x)\big)\sub [-R,R]$ for a.e.\ $x\in \Om$. Further, its barycentre is $A\!:\!\D^2u_\infty$, i.e.
\[
\ \ \ \ \ A :\D^2u_\infty(x)\, =\, \int_\R t\,d [\vartheta_\infty(x) ](t), \ \ \ \text{ a.e. }x\in \Om.
\]
Since $\textbf{H}$ is level-convex, by Jensen's inequality we have
\[
\textbf{H}\left(\int_\R t\,d\big[\vartheta_\infty(x)\big](t)\right)\,\leq\, [\vartheta_\infty(x)]-\underset{t\in \R}{\ess\,\sup}\, \textbf{H}(t)
\]
and since $H$ is continuous and bounded from below
\begin{align}
\ \ \ \ \big\|\textbf{H}\big(A :\D^2u_\infty\big)\big\|_{L^\infty(\Om)}\, &\leq\, \underset{x\in \Om}{\ess\,\sup}\left\{ [\vartheta_\infty(x)]-\underset{t\in \R}{\ess\,\sup}\, \textbf{H}(t)\right\}
\nonumber 
 \end{align}
 
 \begin{align}
\hspace{60pt}  &=\, \lim_{q\ri \infty} \left( \av_\Om  \int_\R \textbf{H}(t)^q \, d\big[\vartheta_\infty(x)\big](t)\,dx \right)^{\frac{1}{q}}
 \label{4.4}
 \\
&\leq\, \underset{q\ri \infty}{\liminf} \left[ \underset{p\ri \infty}{\liminf}\left(   \av_\Om  \int_\R \textbf{H}(t)^q \, d\big[\de_{A:\D^2u_p}(x)\big](t)\,dx   \right)^{\frac{1}{q}}\right]
\nonumber\\
&=\, \underset{q\ri \infty}{\liminf}\left[ \underset{p\ri \infty}{\liminf} \left( \av_\Om   \textbf{H}\big( A:\D^2u_p \big)^q\right)^{\frac{1}{q}}\right].\nonumber
\end{align}
By combining \eqref{4.3}-\eqref{4.4}, the conclusion of (a)-(b) follows. The proof of (c) is identical to that of \cite[Lemma 4]{K8}.  \qed
\ms

Now we establish a complete characterisation of $2$nd order Absolute Minimisers and of the corresponding $\mD$-solutions to the respective equation in the special case of $n=1$.

\begin{theorem}[The fine structure of $2$nd order Absolute Minimisers and of $\mD$-solutions in 1$D$] \label{theorem6} 

Let $\H \in C(\R)$ with $\H \geq \H(0)=0$ and suppose that $\H$ is strictly level convex, that is $\H$ is strictly decreasing on $(-\infty,0)$ and strictly increasing on $(0,\infty)$. Further, suppose that
there exist $C_1,C_2>0$ and $r>1$ such that
\beq
\ \ \ \H(X)\, \geq\, C_1|X|^r \,-\,C_2, \ \ X\in \R.
\eeq
Let $T^-\!<0< T^+$ denote the elements of the level set $\{\H=t\}$ for $t>0$: $\big\{T^-,T^+\big\} = \big\{X\in \R \,:\, \H(X)=t\big\}$. We also suppose that
\beq \label{4.6e}
0\,<\, \underset{t\ri \, \infty}{\lim\inf}\, \left|\frac{\, T^-}{T^+}\right| \, \leq \, \underset{t\ri \, \infty}{\lim\sup} \, \left|\frac{\, T^-}{T^+}\right| \,<\, \infty.
\eeq
Let also $a<b$ in $\R$. We consider the functional \eqref{1.1} on $(a,b)$, that is
\[
\ \ \ \ \E_\infty (u,\mO)\, =\, \big\| \H(u'')\big\|_{L^\infty(\mO)}, \ \ \ u\in W^{2,\infty}(a,b), \ \mO\sub (a,b).
\]
If $\H \in C^1(\R)$ we consider also the corresponding equation \eqref{1.2}, that is
\[
\ \ \ \A^2_\infty u\, =\, \big(\H_X(u'')\big)^3 \big(u'''\big)^2\, =\, 0, \ \ \ \text{ on }(a,b).
\]
Let also $g\in W^{2,\infty}(a,b)$ and set
\beq \label{4.7e}
\mE(g)\, :=\, \frac{g'(b)-g'(a)}{b-a} \, -\, \frac{2\big(g(b)-g(a)-g'(a)(b-a)\big)}{(b-a)^2} .
\eeq
Then, we have:

\smallskip

\noi \emph{(1)} There exists a unique $2$nd order Absolute Minimiser $u_* \in \mathrm{AM}^2(\E_\infty,(a,b))\cap W^{2,\infty}_g(a,b)$ with given boundary values. Moreover, if $\mE(g)=0$ then $u_*$ is a quadratic polynomial function. If $\mE(g)\neq 0$, then $u_*$ is piecewise quadratic with exactly one point $\xi_* \in (a,b)$ at which $u_*''$ does not exist and with $u_*''$ changing sign at $\xi_*$. Further, $\H\big( u_*'' \big)$ extends to a constant function on $(a,b)$. Moreover, $u_*$ coincides with the limit function $u_\infty$ as $p\ri\infty$ of approximate $L^p$ minimisers of Theorem  \ref{theorem4}.

\smallskip

\noi \emph{(2)} Every $2$nd order Absolute Minimiser $u \in \mathrm{AM}^2(\E_\infty,(a,b))$ has the structure described by \emph{(1)} above, i.e.\ $u$ is quadratic if $\mE(u)=0$ and is piecewise quadratic with one point at which $u''$ does not exist and changes sign if  $\mE(u)\neq 0$. Also, $\H(u'')$ extends to a constant function on $(a,b)$ and $u$ coincides with the limit of approximate $L^p$ minimisers in $W^{2,\infty}_u(a,b)$ as $p\ri \infty$.  
\smallskip

\noi \emph{(3)} If in addition $\H\in C^1(\R)$, the Dirichlet problem
\beq \label{4.8e}
\left\{
\begin{array}{rl}
\A^2_\infty u\,=\, 0,\ \ & \text{ in }(a,b), \\
u\, =\, g, \ \ u'\, =\, g',\ \  & \text{ at }\{a,b\}, 
\end{array}
\right.
\eeq
has a unique Absolutely Minimising $\mD$-solution $u_\infty \in W^{2,\infty}_g(a,b)$ which is piecewise quadratic with at most one point in $(a,b)$ at which $u_\infty''$ may not exist. Further, every Absolutely Minimising $\mD$-solution to the problem \eqref{4.8e} is unique has this form.

\smallskip

\noi \emph{(4)} Every $2$nd order Absolute Minimiser of \eqref{1.1} for $n=1$ is a $\mD$-solution to \eqref{1.2}.

\end{theorem}

\begin{remark} \label{remark7} i) The assumption \eqref{4.6e} requires that ``the growth of $\H$ at $+\infty$ can not be too far away from the growth of $\H$ at $-\infty$". It is satisfied for instance if $\H(-X)=\H\big(\al X+o(|X|)\big)$ for some $\al>0$ as $|X|\ri \infty$.

\smallskip

\noi ii) The condition $\mE(g)=0$ (where $\mE(g)$ is given by  \eqref{4.7e}) is necessary and sufficient for the existence of a quadratic polynomial $Q$ with $Q-g \in W^{2,\infty}_0 (a,b)$.  

\smallskip

\noi iii) The converse of item (4) is not true in general (see Section \ref{section8}).
\end{remark}

The proof of Theorem \ref{theorem6} consists of several lemmas.  We begin by recording the following simple observation which relates to Aronsson's result \cite[Lemma 1, p.\ 34]{A1} and is an immediate consequence of the mean value theorem:

\begin{remark} \label{remark8} Suppose that $u$ is a quadratic polynomial on $(\al,\be)\sub \R$ with $u''\equiv C$ and $\phi \in W^{2,\infty}(\al,\be)$ with $\phi \not\equiv u$ and $\phi'=u'$ at $\{\al,\be\}$. Then, there exist measurable sets $A^\pm \sub (\al,\be)$ with $\mL^1(A^\pm)>0$ (positive Lebesgue measure) such that $\phi''$ exists on $A^+ \cup A^-$, whilst we have $\phi''>C$ on $A^+$ and $\phi''<C$ on $A^-$.
\end{remark}

We first consider the much simpler case of $\mE(g)=0$.

\bl \label{lemma9} Every quadratic polynomial $u :\R\ri \R$ is the unique minimiser of $\E_\infty$ over $W^{2,\infty}_u(a,b)$ with respect to its own boundary conditions.
\el

\BPL \ref{lemma9}. Suppose $\phi \in W^{2,\infty}_u(a,b)\set\{u\}$ and $u''\equiv C$ on $(a,b)$. If $C\geq0$, by Remark \ref{remark8} there is a measurable set $A^+\sub (a,b)$ with $\mL^1(A^+)>0$ such that $\phi''>C\geq 0$ on $A^+$. Since $\H$ is strictly increasing on $(0,\infty)$, we have 
\[
\E_\infty\big(\phi,(a,b)\big)\, \geq\, \E_\infty\big(\phi,A^+\big)\, =\, \big\|\H(\phi'')\big\|_{L^\infty(A^+)} >\, \H(C)\, =\, \E_\infty\big(u,(a,b)\big).
\]
The case $C<0$ follows analogously since by Remark \ref{remark8} there is a measurable $A^-\sub (a,b)$ with $\mL^1(A^-)>0$ and $\phi''<C<0$, whilst $\H$ is strictly decreasing on $(-\infty,0)$ and we again obtain $\E_\infty\big(\phi,(a,b)\big) > \E_\infty\big(u,(a,b)\big)$. Uniqueness follows by the strictness of the energy inequalities. \qed

\ms

We now consider the case of $\mE(g)\neq 0$.

\bl \label{lemma10} If $\mE(g)\neq 0$, there exists a unique piecewise quadratic function $u_*\in W^{2,\infty}_g(a,b)$ and a $\xi_* \in (a,b)$ such that $u''\equiv L_*$ on $(a,\xi_*)$, $u''\equiv R_*$ on $(\xi_*,b)$ and 
\[
\ \ \ \big\| \H(u_*'')\big\|_{L^\infty(a,b)} =\, \max\big\{\H(L_*),\H(R_*)\big\}, \ \ \ \H(L_*)\,=\,\H(R_*), \ \ L_*R_*<\,0.
\]
\el

\BPL \ref{lemma10}. For brevity we set $g(a)=A$, $g(b)=B$, $g'(a)=A'$ and $g'(b)=B'$. We first show that a piecewise quadratic function in $W^{2,\infty}_g(a,b)$ with one matching point indeed exists. We fix parameters $R,L \in \R$ and set
\[
\left\{ \ \ 
\begin{split}
u^L(x)\,&:=\, A \,+\,A'(x-a)\,+\, \frac{L}{2}(x-a)^2, \\
u^R(x)\,&:=\, B \,+\,B'(x-b)\,+\, \frac{R}{2}(x-b)^2.
\ \ \ \ 
\end{split}
\right.
\]
The condition of $C^1$ matching of $u^L,u^R$ at some ${\xi} \in (a,b)$ to a single function $u=u^L\chi_{(a,{\xi})}+ u^R\chi_{[{\xi},b)}$ in $W^{2,\infty}_g(a,b)$ is equivalent to
\begin{align}
A \,+\,A'({\xi}-a)\,+\, \frac{L}{2}({\xi}-a)^2\, &=\, B \,+\,B'({\xi}-b)\,+\, \frac{R}{2}({\xi}-b)^2, \label{4.10e} \\
(L-R){\xi}\, &=\, B'\,-\,A' \,+\, aL\,-\,bR .    \label{4.11e}
\end{align}
Since by assumption $\mE(g)\neq 0$, it follows that $R\neq L$. We now set
\beq \label{4.12e}
C_0\,:=\, \frac{B'-A'}{b-a},\ \ \ C_1\,:=\, \frac{2\big(A-B-B'(a-b)\big)}{(b-a)^2}, \ \ \ C_2\,:=\, \frac{2\big(B-A-A'(b-a)\big)}{(b-a)^2}.
\eeq
By \eqref{4.11e}-\eqref{4.12e}, the constraint $a<{\xi}<b$ is equivalent to
\beq \label{4.13e}
\left\{ \ \ \
\begin{split}
& R\,<\, C_0\,<\, L, \ \ \ \text{ if }R<L,\\ 
& L\,<\, C_0\,<\, R, \ \ \ \text{ if }L<R.
\end{split}
\right.
\eeq
By cancelling $\xi$ from \eqref{4.10e} with the aid of \eqref{4.11e}, we obtain that the admissible pairs $(R,L)\in \R^2$ for matching lie on the hyperbola $\mC \sub \R^2$ given by the equation
\beq \label{4.14e}
\big(R\,-\,C_1\big)\, \big(L\,-\,C_2\big)\, =\, C_1C_2\,-\,C_0^2.
\eeq
Then, \eqref{4.14e} coincides with the condition of vanishing discriminant of the algebraic equation \eqref{4.10e} (when considered as a binomial equation with respect to ${\xi}$) and the unique point ${\xi}$ for matching is that given by \eqref{4.11e}. Note that the lines $\{R=C_1\}$ and $\{L=C_2\}$ are the asymptotes of $\mC$. Since
\[
\mE(g)\, =\, \frac{B'-A'}{b-a}\, -\, \frac{2\big(B-A-A'(b-a)\big)}{(b-a)^2}\,=\, C_0-C_2
\]
the following facts can be easily verified by elementary algebraic calculations:
\begin{align}
\mE(g)>\, \,0  \ \ \   \Longleftrightarrow & \ \ \ C_1\,>\,C_0\,>\,C_2, \label{4.15e} \\
C_0^2 \, -\, C_1C_2\, & =\,\mE(g)^2 \,>0,  \label{4.16e} \\
\ \ (C_0,C_0) \, \in \mC \ \& \ \{R=L\} & \text{ is tangent to $\mC$ at $(C_0,C_0)$.}  \label{4.17e}
\end{align}
Note further that $\mC$ lies on the 2nd and 4th quadrants of $\R^2$, i.e. $\big(R-C_1\big)\big(L-C_2\big)<0$ for all $(R,L) \in \mC$. We now derive the remaining constraints that $(R,L)$ have to satisfy in order to be admissible. We set $\de:=b-{\xi}$ and rewrite \eqref{4.10e}-\eqref{4.11e} as
\[
\left\{\ 
\begin{split}
\de^2(L-R)\, -\, 2\de\Big[ L(b-a)-(B'\,-\,A') \Big] + L(b-a)^2\, &=\, 2\big(B-A-A'(b-a)\big), \\
 \de (L-R)\, &=\, L(b-a)-(B'\,-\,A').    
\end{split}
\right. \ \ 
\]
The above imply the inequality
\beq \label{4.18e}
\de^2 \, =\, \frac{(b-a)^2}{L-R}\left(L\,-\, \frac{2\big(B-A-A'(b-a)\big)}{(b-a)^2}\right)\,\geq\, 0.
\eeq
Similarly, we set $\e:={\xi}-a$ and rewrite \eqref{4.10e}-\eqref{4.11e} as
\[
\left\{
\begin{split}
\e^2(L-R)\, +\, 2\e\Big[ R(b-a)-(B'\,-\,A') \Big]-R(b-a)^2 &= -2\big(A-B-B'(a-b)\big), \\
 -\e (L-R)\, &=\, R(b-a)-(B'\,-\,A').    
\end{split}
\right. \ \ 
\]
The above imply the inequality
\beq \label{4.19e}
\e^2 \, =\, \frac{(b-a)^2}{R-L}\left(R\,-\, \frac{2\big(A-B-B'(a-b)\big)}{(b-a)^2} \right)\,\geq\, 0.
\eeq
By \eqref{4.19e}, \eqref{4.18e}, \eqref{4.13e} and \eqref{4.12e}, we have that the admissible pairs $(R,L)$ lie on the constraint set
\beq \label{4.20e}
\mK\,:=\, \bigg\{ (R,L)\in \R^2\, \bigg| \, 
\begin{array}{l}
R\,\leq\, C_1,\ \,L\,\geq\, C_2 \,\text{ and } R<C_0<L, \ \text{ if } R<L;\\
R\,\geq\, C_1, \ \, L\,\leq\, C_2  \,\text{ and } L<C_0<R, \ \text{ if }\, L<R 
\end{array}
\bigg\} . 
\eeq
\[
\underset{\text{Figure 1.}}{\includegraphics[scale=0.21]{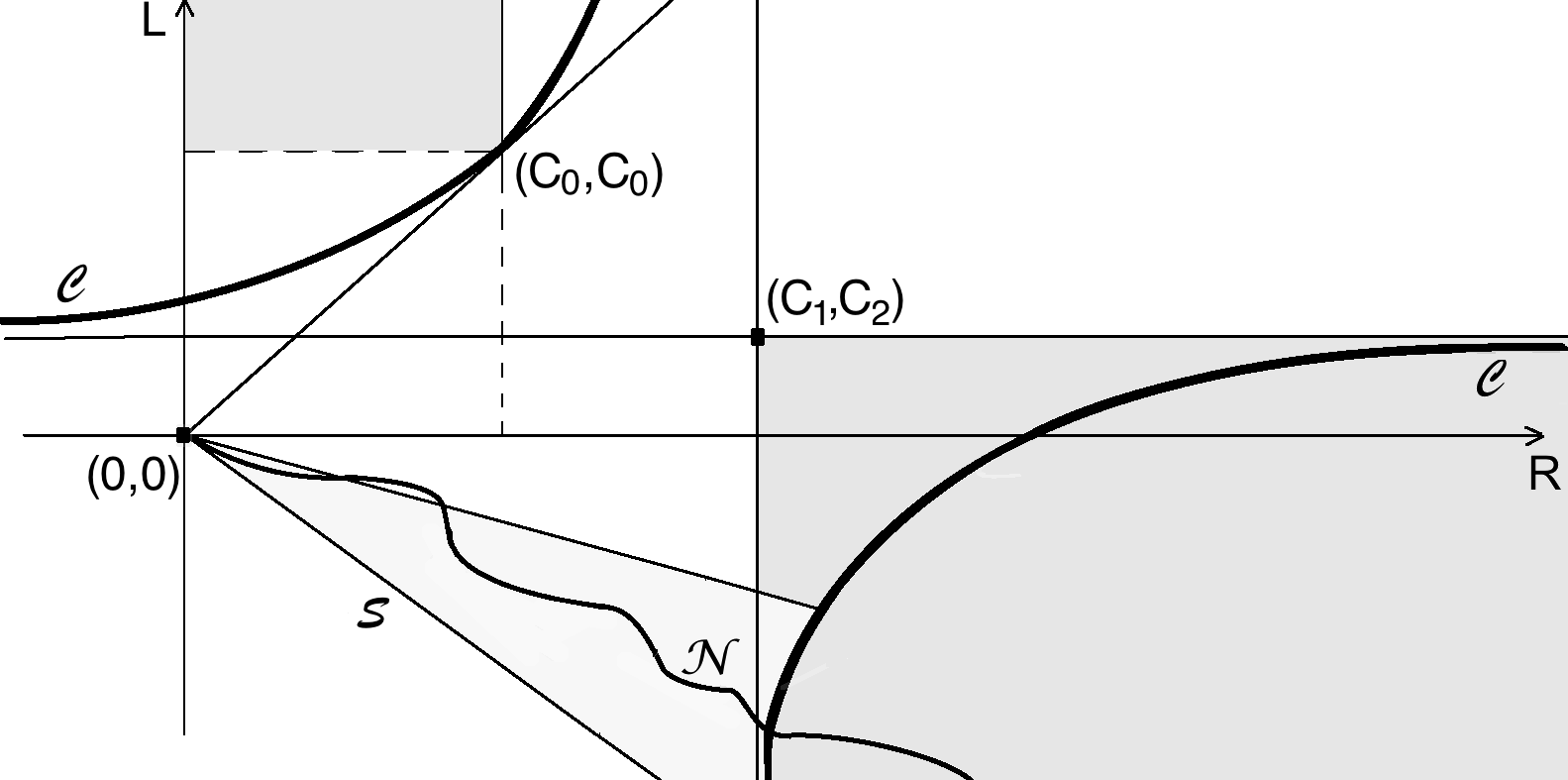}}
\]
We consider first the case $\mE(g)>0$. By \eqref{4.15e}, the set of admissible pairs is
\[
\mC \cap \mK\, =\, \Big\{(R,L)\in \mC\, :\, R\geq C_1,\ L\leq C_2 \Big\}.
\]
Consider now the continuous curve
\[
\mN\,:=\, \Big\{ (R,L)\in \R^2 \,:\, L<0<R, \, \H(L)=\H(R)\Big\}.
\]
By our assumptions on $\H$, the correspondences $R\mapsto L$ and $L \mapsto R$ are inverse of each other and hence $\mN$ is the graph of a monotone function. By our assumption \eqref{4.6e}, for $|(R,L)|$ large enough $\mN$ lies in a sector of the form
\[
\S\,:=\, \Big\{(R,L)\in \R^2\ : \ \Big|\angle \Big((R,L),(1,-1) \Big) \Big|\, \leq \, a  \Big\}
\]
for some $a<\pi/4$. Since $\mC \cap \mK$ intersects every ray in the sector $\S$ emanating from the origin, it follows that there is a unique point of intersection $(R_*,L_*) \in \mN\cap \mC \cap \mK$ giving rise to a unique matching point $\xi_*$ given by \eqref{4.11e}. Then, the function $u_*:=\, \smash{u^{L_*}}\chi_{(a,\xi_*)}+ \smash{u^{R_*}}\chi_{[\xi_*,b)}$ has the desired properties. The case of $\mE(g)<0$ follows analogously, so the proof of the lemma is complete.   \qed

\bl \label{lemma11} If $\mE(g)\neq 0$, the function $u_*\in W^{2,\infty}_g(a,b)$ obtained in Lemma \ref{lemma10} is the unique minimiser of $\E_\infty$ over $W^{2,\infty}_g(a,b)$.
\el

\BPL \ref{lemma11}. We simplify the notation and symbolise $u_*$ by just $u$. Fix $v \in W^{2,\infty}_g(a,b) \set\{u\}$. Then, $v' \in W^{1,\infty}_{u'}(a,b)$ and $v' \not\equiv u'$. Without loss of generality consider the case $R<0<L$. If $v'$ intersects $u'$ at some point in $(a,\xi)$, then by Remark \ref{remark8} there exists $A^+\sub (a,\xi)$ with $\mL^1(A^+)>0$ such that $v'' >L>0$ on $A^+$. Since $\H$ is strictly increasing on $(0,\infty)$, we have
\[
\begin{split}
\E_\infty\big(v,(a,b) \big)\, &\geq\, \E_\infty\big(v,A^+ \big) 
\,
=\, \underset{A^+}{\ess\,\sup}\, \H(v'')\, >
\\
&>\, \H(L)
\, =\, \max\big\{\H(L),\H(R) \big\}
\,=\, \E_\infty\big(u,(a,b) \big).
\end{split}
\]
Similarly, if $v'$ intersects $u'$ at some point in $[{\xi},b)$, then by Remark \ref{remark8} there exists $A^-\sub ({\xi},b)$ with $\mL^1(A^-)>0$ such that $v'' <R<0$ on $A^-$. Since $\H$ is strictly decreasing on $(-\infty,0)$, we again obtain $\E_\infty\big(v,(a,b) \big) > \E_\infty\big(u,(a,b) \big)$. Hence, it remains to consider the cases that $v'$ lies either above or below $u'$ over $(a,b)$. The former case can also be handled by Remark \ref{remark8}: Since $v'(a)=u'(a)$ and $v'({\xi})>u'({\xi})$, there is a set $A^+\sub (a,{\xi})$ with $\mL^1(A^+)>0$ such that
\[
\ \ \ v''\, >\, \frac{v'({\xi})-v'(a)}{{\xi} -a }\, >\, \frac{u'({\xi})-u'(a)}{{\xi} -a }\, =\, L\,>\,0, \ \ \ \text{ on }A^+,
\] 
and again $\E_\infty\big(v,(a,b) \big) > \E_\infty\big(u,(a,b) \big)$. Finally, if $v'<u'$ on $(a,b)$, then since $u(a)=v(a)$ by integration we get $v<u$ on $(a,b)$. Since $u$ is quadratic on $({\xi},b)$ with $u''=R$ and $v({\xi})<u({\xi})$, by Taylor's theorem we have
\[
\begin{split}
v(b)\, +\, v'(b)({\xi}-b)\,+\ & \frac{1}{2}\left( \int_0^1 v''\big(b+ t ({\xi}-b) \big)\big[2(1-t) \big] \, dt \right)({\xi}-b)^2 \\
& < \, u(b)\, +\, u'(b) ({\xi}-b)\, +\, \frac{R}{2}({\xi}-b)^2.
\end{split}
\]
Since $u(b)=v(b)=B$, $u'(b)=v'(b)=B'$, by considering the absolutely continuous probability measure $\mu<<\mL^1$ on $[0,1]$ given by $\mu(E):=\int_E 2(1-t) d t$, we deduce
\[
\ \  \int_0^1 v''\big( b+ t ({\xi}-b)  \big) \, d\mu(t) \, <\, R \, \mu\big([0,1]\big).
\]
Hence, there exists a measurable set $A\sub (0,1)$ with $\mL^1(A)>0$ such that $v''\big(b+ t ({\xi}-b) ) \big)<R$ for all points $t \in A$. Since $R<0$, by arguing as before we obtain $\E_\infty\big(v,(a,b) \big) > \E_\infty\big(u,(a,b) \big)$. The lemma ensues.   \qed

\bl \label{lemma12} In the setting of Theorem \ref{theorem4}, the function $u_\infty$ constructed therein is a $2$nd order absolute minimiser of $\E_\infty$ in $W^{2,\infty}_g(a,b)$.
\el

\BPL \ref{lemma12}. We begin with the following observation: given $A,B,A',B',$ $a,b \in \R$ with $a<b$, the unique cubic Hermite interpolant $Q:\R\larrow \R$ satisfying
$Q(a)=A$, $Q(b)=B$, $Q'(a)=A'$ and $Q'(b)=B'$. Further, let $(v_p)_1^\infty \sub W^{2,\infty}(a,b)$ be any sequence of functions satisfying $v_p\larrow v_\infty$ in $C^1[a,b]$ as $p\ri \infty$. If $Q_p$ is the cubic polynomial such that $Q_p-v_p \in W^{2,\infty}_0(a,b)$ for $p\in \N\cup\{\infty\}$, namely
\[
\text{$Q_p(a)\,=\,v_p(a)$, \ \, $Q_p(b)\,=\,v_p(b)$, \ \, $Q_p'(a)\,=\,v_p'(a)$, \ \, $Q_p'(b)\,=\,v_p'(b)$,}
\]
we have $Q_p \larrow Q_\infty$ \texttt{strongly} in $W^{2,\infty}(a,b)$ as $p\ri \infty$.

\smallskip

We now continue with the existence of an absolute minimiser. Let $(u_p)_{p=1}^\infty$ be the sequence of approximate minimisers of Theorem \ref{theorem4} which satisfies $u_p \weak u_\infty$ in $W^{2,q}(\Om)$ as $p\ri \infty$ along a subsequence for any $q>1$. Fix an $\Om'\Subset \Om$ and $\phi \in W^{2,\infty}_0(\Om')$. Since any open set on $\R$ is a countable disjoint union of intervals, we may assume $\Om'=(a,b)$. In order to conclude, it suffices to show that
\beq \label{4.6}
\ \ \ \E_\infty\big(u_\infty,(a,b) \big)\,\leq\, \E_\infty\big(u_\infty +\,\phi,(a,b) \big),
\eeq
when $\phi \in W^{2,\infty}_0(a,b)$. Consider for any $p\in \N\cup\{\infty\}$ the unique cubic polynomial such that $Q_p-u_p \in W^{2,\infty}_0(a,b)$. By the above observations and Theorem \ref{theorem4}(a)-(b), along a subsequence $p_j \ri \infty$ we have 
\beq \label{4.7}
\ \ \ \ \ Q_p\larrow Q_\infty \ \  \text{ in }W^{2,\infty}(a,b),\ \ \text{ as }p\ri \infty.
\eeq
We define for $p\in\N$ the function $\phi_p := \phi \, +\, u_\infty -\, u_p\, +\, Q_p -\, Q_\infty$. Since all three functions $\phi$, $u_\infty-Q_\infty$ and  $u_p-Q_p$ are in $W^{2,\infty}_0(a,b)$, the same is true for $\phi_p$ as well. By Theorem \ref{theorem4}(b) and by the additivity of the integral, we have
\beq \label{4.8}
\begin{split}
\E_p\big(u_p,(a,b)\big)\, \leq & \ \, 2^{-p}\, +\, \E_p\big(u_p+\phi_p\, ,(a,b)\big)
\\
= & \ \, 2^{-p}\, +\,  \E_p\Big(u_\infty+\phi+[Q_p-Q_\infty]\, , (a,b)\Big)
\\
\leq & \ \, 2^{-p} \, +\,  \left(\frac{b-a}{\mL^1(\Om)}\right)^{\!1/p}\E_\infty\Big(u_\infty+\phi+[Q_p-Q_\infty]\, , (a,b)\Big).
\end{split}
\eeq
By invoking \eqref{4.7} and passing to the limit in \eqref{4.8} as $p\ri \infty$, we deduce  \eqref{4.6} as a consequence of Theorem \ref{theorem4}(c). The lemma ensues.          \qed
\ms

We can now prove the result by using the above lemmas.

\BPT \ref{theorem6}. (1) By Lemmas \ref{lemma9}-\ref{lemma11}, there exists a unique minimiser $u_*$ of $\E_\infty$ in $W^{2,\infty}_g(a,b)$ which is either quadratic or piecewise quadratic with at most one breaking point for $u_*''$ at which it changes sign and $\H\big(u_*''\big)$ extends to a constant function on $(a,b)$. By Lemma \ref{lemma12}, the limit $u_\infty$ as $p\ri \infty$ of approximate $L^p$ minimisers is a $2$nd order Absolute Minimiser and a fortiori a minimiser of $\E_\infty$ in $W^{2,\infty}_g(a,b)$. Thus, $u_*\equiv u_\infty$ and this is the unique element of $\mathrm{AM}^2\big(\E_\infty,(a,b)\big) \cap W^{2,\infty}_g(a,b)$. (2) is a consequence of (1). (3) and (4) follow by the fact that $\H\big( u_*'' \big)=C$ a.e.\ on $(a,b)$, uniqueness and Claim \ref{claim} of Section \ref{section7} that follows. The theorem ensues. \qed

\section{Variational characterisation of $\A^2_\infty $ via $\mathrm{AM}^2$ \& uniqueness}  \label{section6}

Herein we show that $2$nd order Absolute Minimisers of \eqref{1.1} in $C^3(\Om)$ solve the fully nonlinear PDE \eqref{1.2}. The converse is also true in the case that $\H$ depends on the hessian via a scalar projection of it along a matrix. Let us note that in that case by the Spectral Theorem and a change of variables, the study of this functional can be reduced to the study of one depending on the hessian via \emph{just} the Laplacian, but we find it more elucidating to retain this seemingly more general form. As a consequence, we deduce the uniqueness of $C^3$ $2$nd order Absolute Minimisers and of classical solutions to the PDE. Accordingly, the main result of this section is:

\begin{theorem}[Variational characterisation of $\A^2_\infty $ via $\E_\infty$] \label{theorem5}

Given a non-negative function $\H\in C^1\big(\R^{n^{\ot 2}}_s\big)$, an open set $\Om\sub \R^n$  and $u\in C^3(\Om)$, consider the supremal functional $\E_\infty$ given by \eqref{1.1} and the fully nonlinear equation $\A^2_\infty  u=0$ given by \eqref{1.2}. Then:

\ms

\noi \emph{(A)} If $u \in \mathrm{AM}^2(\E_\infty,\Om) \cap C^3(\Om)$, namely if it is a $C^3$ $2$nd order Absolute Minimiser (Definition \ref{definition1}), then it solves $\A^2_\infty  u=0$ on $\Om$.

\ms

\noi \emph{(B)} Suppose that $\Om$ is connected and $\H$ has the form $\H(X)\, =\, \textbf{\emph{H}}\big(A\!:\!X\big)$ for a fixed (strictly) positive matrix $A \in \R^{n^{\ot 2}}_s$ and some level-convex function $\textbf{\emph{H}} \in C^1(\R)$ (that is $\textbf{\emph{H}}$ has for any $t\geq 0$ convex sub-level sets $\{\textbf{\emph{H}}\leq t\}$) such that $\{\textbf{\emph{H}}=t\}$ consists of at most $2$ points.  Then, the statements \emph{(1)-(4)} below are equivalent:
\begin{enumerate}
\item $u \in \mathrm{AM}^2(\E_\infty,\Om)$.
\smallskip

\item There exists $C\geq 0$ such that $ A:\D^2u \equiv C$ in $\Om$.

\smallskip

\item There exists $c\geq 0$ such that $\textbf{\emph{H}}\big(A:\D^2u\big)\equiv c$ in $\Om$.

\item $\A^2_\infty u =0$ in $\Om$.
\end{enumerate}

\end{theorem}

As a consequence, we deduce the next result:

\begin{corollary}[Uniqueness of $C^3$ Absolute Minimisers and of $C^3$ solutions to the Dirichlet problem] \label{corollary6}

In the setting of Theorem \ref{theorem5}(B) above, suppose in addition $\Om$ is bounded. Then, for any $g\in W^{2,\infty}(\Om)$, the intersection
\[
C^3(\Om)\cap W^{2,\infty}_g(\Om) \cap \mathrm{AM}^2(\E_\infty,\Om)
\]
contains at most one element, namely there is at most one $C^3$ Absolute Minimiser $u$ of $2$nd order which satisfies $u=g$ and $\D u=\D g$ on $\p\Om$. Further, the problem
\beq
\left\{
\begin{array}{rl}
\Big(\textbf{\emph{H}}'\big(A:\D^2u\big)\, A\big)^{\ot 3}: \big(\D^3u\big)^{\ot 2}\, =\, 0,\ \ & \ \text{ in }\Om,
\\
u\,=\,g, \ \D u\,= \, \D g,\  \ & \ \text{ on }\p\Om,\ms
\end{array}
\right.
\eeq
has at most one solution in $C^3(\Om)\cap W^{2,\infty}_g(\Om)$.  
\end{corollary}

We begin with a simple lemma which is relevant to some of the results of \cite{K10}.

\begin{lemma} \label{lemma8} Let $\H\in C^1\big(\R_s^{n^{\ot 2}}\big)$ be given and consider the functional \eqref{1.1}. Let also $\Om\sub \R^n$ be an open set.

\smallskip

\noi \emph{(a)} For any $u\in C^2(\Om)$ and $\mO \Subset \Om$, we set 
\beq
\mO(u)\,:=\, \bigg\{x\in \overline{\mO}\ :\ \H\big(\D^2u(x)\big)\, =\, \big\|\H\big(\D^2u\big)\big\|_{L^\infty(\mO)} \bigg\}.
\eeq
Then, if $u \in \mathrm{AM}^2(\E_\infty,\Om)$ (Definition \ref{definition1}), it follows that
\[
\ \ \ \H_X\big(\D^2u\big) :\D^2\phi\, =\, 0, \ \ \text{ on }\mO(u),
\]
for any $\phi \in W^{2,\infty}_0(\mO)\cap C^2(\overline{\mO})$.

\smallskip

\noi  \emph{(b)} For any $u\in C^2(\Om)$, $x\in \Om$ and $0<\e<\frac{1}{2}\dist(x,\p\Om)$, we set \beq \label{6.3}
\Om_\e(x)\,:=\, \mB_\e(x)\cap \Big\{\H\big(\D^2u\big) < \H\big(\D^2u(x)\big)\Big\}
\eeq
and suppose $\Om_\e(x)\neq \emptyset$. Then, $x\in \p(\Om_\e(x))$ and
\[
\H\big(\D^2u(y)\big)\, =\, \big\|\H\big(\D^2u\big)\big\|_{L^\infty(\Om_\e(x))}
\]
for all $y \in \mB_\e(x) \cap \p(\Om_\e(x))$ (in particular for $y=x$).

\end{lemma}

\BPL \ref{lemma8}. (a) This is an application of Danskin's theorem \cite{Da}. Fix $\mO\Subset \Om$ and $\phi \in W^{2,\infty}_0(\mO)\cap C^2(\overline{\mO})$. If $u \in \mathrm{AM}^2(\E_\infty,\Om)$, then the real function $t\mapsto \|\H\big(\D^2u +t\D^2\phi\big)\|_{L^\infty(\mO)}$ attains a local minimum at $t=0$, so if its derivative exists it must vanish. By Danskin's theorem we have
\[
\begin{split}
\frac{d}{dt}\Big|_{t=0} & \bigg(\big\|\H\big(\D^2u +t\D^2\phi\big)\big\|_{L^\infty(\mO)}\bigg)\, =\, \frac{d}{dt}\Big|_{t=0}\bigg(\max_{\overline{\mO}}\, \H\big(\D^2u +t\D^2\phi\big) \bigg)\, =
\\
&=\, \max_{x\in  {\mO(u)}}\bigg(\frac{d}{dt}\Big|_{t=0} \H\big(\D^2u +t\D^2\phi\big)(x)\bigg)
\, =\, \max_{x\in  {\mO(u)}} \bigg( \H_X\big(\D^2u(x)\big):\D^2\phi(x) \bigg)
\end{split}
\]
and upon replacing $\phi$ with $-\phi$, the conclusion follows.  

\smallskip

\noi (b) This is obvious from the definitions. \qed

\ms

\BPT \ref{theorem5}. (A) Fix $x\in \Om$ and $\e\in(0,\dist(x,\p\Om))$. If $\Om_\e(x)= \emptyset$, then we have $\H(\D^2u)\geq \H(\D^2u(x))$ on $\mB_\e(x)$, which implies that $x$ is an interior minimum, therefore giving $\D\big(\H(\D^2u))(x)=0$. Hence, $\A^2_\infty u(x)=0$. Consequently, we may assume $\Om_\e(x)\neq \emptyset$. By (a) and (b) of Lemma \ref{lemma8} above, if $u \in \mathrm{AM}^2(\E_\infty,\Om) \cap C^2(\Om)$, we have  
\beq \label{6.4}
\H_X\big(\D^2u(x)\big) :\D^2\phi(x)\, =\, 0
\eeq
for any $\phi \in W^{2,\infty}_0(\Om_\e(x))\cap C^2(\overline{\Om_\e(x)})$. We illustrate the idea of the rest of the proof by showing it first under one more degree of regularity.

\begin{claim} \label{claim9} Theorem \ref{theorem5}(A) is true if in addition $u \in C^4(\Om)$ and $\H \in C^2\big(\R^{n^{\ot 3}}_s\big)$.
\end{claim}

\BPC \ref{claim9}. Fix $x,\e,\phi$ as above. If $\Om_\e(x)\neq \emptyset$, let $\ze \in C^\infty_c(\mB_\e(x))$ be a cut off function with $\ze\geq 0$ on $\mB_{\e}(x)$ and $\ze\equiv 1$ on $\mB_{\e/2}(x)$. We define
\beq \label{6.5}
\phi \,:=\, \frac{1}{2} \Big(\H\big(\D^2u \big)-\H\big(\D^2u(x)\big)\Big)^{\!2} \zeta .
\eeq
Then, $\phi \in C^2(\overline{\Om_\e(x)})$, $\phi =0$ on $\p(\Om_\e(x))$ and on $\overline{\Om_\e(x)}$ we have
\beq \label{6.6}
\begin{split}
\D\phi = \ze \Big(\H\big(\D^2u \big)-\H\big(\D^2u(x)\big)\Big) \D\big(\H\big(\D^2u \big)\big) + \Big(\H\big(\D^2u \big)-\H\big(\D^2u(x)\big)\Big)^{\!2}\,\D\zeta
\end{split}
\eeq
whilst on $\mB_{\e/2}(x)$ we have
\[
\D^2\phi \,=\, \D\big(\H\big(\D^2u \big)\big)\ot \D\big(\H\big(\D^2u \big)\big)\,+\, \Big(\H\big(\D^2u \big)-\H\big(\D^2u(x)\big)\Big) \D^2\big(\H\big(\D^2u \big)\big).
\]
Hence, $\D\phi =0$ on $\p(\Om_\e(x))$ and 
\[
\D^2\phi(x) \,=\, \D\big(\H\big(\D^2u  \big)\big)(x) \ot \D\big(\H\big(\D^2u  \big)\big)(x).
\]
By inserting the above to \eqref{6.4}, we obtain
\[
\ \ \A^2_\infty u(x)\,=\,\H_X\big(\D^2u(x)\big) : \D\big(\H\big(\D^2u  \big)\big)(x) \ot \D\big(\H\big(\D^2u  \big)\big)(x)\, =\,0,
\]
for any $x \in \Om$. The claim ensues. \qed
\ms

Now we complete the proof of (A) by considering the case of ``natural" regularity $u \in \mathrm{AM}^2(\E_\infty,\Om) \cap C^3(\Om)$ and $\H \in C^1\big(\R^{n^{\ot 3}}_s\big)$. Let again $\phi$ and $\ze$ be as in \eqref{6.5}. Then, we have $\phi \in C^1(\overline{\Om_\e(x)})$, $\phi =0$ on $\p(\Om_\e(x))$ and \eqref{6.6} still holds giving $\D\phi =0$ on $\p(\Om_\e(x))$. The problem is that $ \D^2\big(\H\big(\D^2u \big)\big)$ may not exist. 

However, we now show that \emph{$\H\big(\D^2u \big)$ is twice differentiable in the sense of Whitney} (see \cite{W, M, F}) on the closed set $\overline{\mB}_\e(x)\cap \p(\Om_\e(x))$ and hence it admits a $C^2$ extension which coincides with it up to $2$nd order on $\overline{\mB}_\e(x)\cap \p(\Om_\e(x))$. To this end we need to introduce some notation. Given any functions $f,g :\R^n \larrow \R$, we set
\[
\ \ \ \ \ \ \D^{1,z}f(y)\,:=\, \frac{f(y+z)-f(y)}{|z|}, \ \ \ f \, :\ \R^n \larrow \R,\ z\neq 0,
\]
and note the elementary identity
\beq \label{6.7a}
\ \ D^{1,z}(fg)(y)\, =\, f(y)\,\D^{1,z}g(y)\, +\, g(y+z)\,\D^{1,z}f(y).
\eeq
Let us also set for brevity $\mathscr{H}:= \H(\D^2u)\,-\,\H\big(\D^2u(x)\big)$ and $\mK:=\overline{\mB}_\e(x)\cap \p(\Om_\e(x))$. By applying \eqref{6.7a}, we have
\[
\begin{split}
\D^{1,z}\D\Big(\frac{1}{2}\mathscr{H}^2\Big)\, &= \D^{1,z}\big(\mathscr{H}\,\D\mathscr{H})
\,= \, \mathscr{H} \, \big(\D^{1,z}\D\mathscr{H}\big)\, +\, \big(\D^{1,z}\mathscr{H}\big)\, \D\mathscr{H}(\cdot +z) .
\end{split}
\]
Note now that \eqref{6.3} yields $\mathscr{H}\equiv 0$ on $\mK$. Hence, we obtain
\beq \label{6.8a}
\ \ \ \D^{1,z}\D\Big(\frac{1}{2}\mathscr{H}^2\Big)\, \,= \big(\D^{1,z}\mathscr{H}\big)\, \D\mathscr{H}(\cdot +z), \ \ \text{ on }\mK.
\eeq
We now claim that there exists an increasing modulus of continuity $\om\in C[0,\infty)$ with $\om(0)=0$ such that
\beq \label{6.9a}
\max_{y\in\mK}\, \left|\frac{\D (\frac{1}{2}\mathscr{H}^2)(z+y)-\D(\frac{1}{2}\mathscr{H}^2)(y)-(\D \mathscr{H}(y)\ot \D \mathscr{H}(y)\big)z}{|z|} \right| \, \leq \,\om(|z|) 
\eeq
for any $0<|z|<\e/2$. In order to establish \eqref{6.9a}, we fix a $y\in \mK$ and calculate:
\[
\begin{split}
&\bigg| \frac{\D (\mathscr{H}^2)(z+y) -\D(\mathscr{H}^2)(y)-2(\D \mathscr{H}(y)\ot \D \mathscr{H}(y)\big)z}{|z|} \bigg|
 \\
 & \hspace{30pt}=\ \left|\D^{1,z}\D(\mathscr{H}^2)(y)-2\frac{\D \mathscr{H}(y) \ot \D \mathscr{H}(y)}{|z|}\,z \right|
 \\
 &  \hspace{25pt}\overset{\eqref{6.8a}}{=}\,  2\left|  \D^{1,z}\mathscr{H}(y)\, \D\mathscr{H}(y +z)-\Big( \D \mathscr{H}(y) \ot \D \mathscr{H}(y) \Big)\frac{z}{|z|} \right|
 \\
 &  \hspace{30pt}=\  2\left|  \frac{ \mathscr{H}(y+z)- \mathscr{H}(y)}{|z|}\, \D\mathscr{H}(y +z)-\left(\D \mathscr{H}(y)\cdot\frac{z}{|z|}\right)\D \mathscr{H}(y) \right|
\\
&  \hspace{30pt}\leq \  2\big|\D \mathscr{H}(y) \big| \Big|\D \mathscr{H}(y+z)-\D\mathscr{H}(y) \Big| 
\\
&  \hspace{40pt} +\ 2 \big|\D \mathscr{H}(y+z) \big| 
\left| \frac{ \mathscr{H}(y+z) -\mathscr{H}(y)- \D\mathscr{H}(y)\cdot z }{|z|} \right|
\\
&  \hspace{30pt}\leq \ \big\|\D \mathscr{H}\big\|_{L^\infty(\mB_{2\e}(x))}\Big(\om_1(|z|)\,+\, \om_2(|z|)\Big),
\end{split}
\]
for some moduli of continuity $\om_1,\om_2$ (by the $C^1$-regularity of $\mathscr{H}$).
Hence, \eqref{6.9a} has been established. Further, since $\D(\mathscr{H}^2)\equiv 0$ on $\mK$ we obviously have
\beq \label{6.10a}
\max_{y\in \mK}\,\Big|\mathscr{H}^2(y+z) -\mathscr{H}^2(y) \Big|\, =\, o(|z|)
\eeq
as $z\ri 0$, while also
\beq \label{6.11a}
\ \ \max_{y\in\mK}\, \left| \mathscr{H}^2(z+y)-\mathscr{H}^2(y)-\Big(\D \mathscr{H}(y)\ot \D \mathscr{H}(y)\Big):z\ot z  \right| \, \leq \, o(|z|^2),
\eeq
as $z\ri 0$. The inequality \eqref{6.11a} is an easy consequence of \eqref{6.9a} and the identity
\[
\mathscr{H}^2(w+y)-\mathscr{H}^2(y)\, =\, w \cdot \int_0^1 \D(\mathscr{H}^2)\big(y +\la w\big)\,d\la .
\]
Conclusively, by \eqref{6.9a}-\eqref{6.11a} the function $\mathscr{H}^2$ is twice (Whitney) differentiable on the closed set $\mK$ with Whitney hessian
\[
\D^2\Big(\frac{1}{2}\mathscr{H}^2\Big) \,=\, \D\big(\H\big(\D^2u \big)\big)\ot \D\big(\H\big(\D^2u\big)\big).
\]
By the Whitney extension theorem (\cite{W, M, F}), there exists an extension $\Phi 
\in C^2(\R^n)$ such that, on $\mK$ we have $\Phi = \mathscr{H}^2/2 \equiv 0$, $\D\Phi = \D(\mathscr{H}^2/2) \equiv 0$ and
\[
\D^2\Phi \, =\, \D^2 \big(\mathscr{H}^2/2\big)=\D\big(\H\big(\D^2u \big)\big)\ot \D\big(\H\big(\D^2u\big)\big).
\]
Let now $\ze$ be the cut-off function of Claim \ref{claim9}. The test function $\phi:=\Phi \ze$ satisfies $\phi \in W_0^{2,\infty}(\Om_\e(x))\cap C^2(\overline{\Om_\e(x)})$ and 
\[
\D^2\phi(x)\, =\, \D\big(\H\big(\D^2u \big)\big)(x)\ot \D\big(\H\big(\D^2u\big)\big)(x).
\]
By inserting the above to \eqref{6.4} we conclude that $\A^2_\infty u(x)=0$ for the arbitrary point $x\in \Om$.

\ms

\noi (B)  (4) $\Rightarrow$ (3): We rewrite the equation $\A^2_\infty u =0$ as
\[
\textbf{H}'\big(A:\D^2u\big) \Big[A: \D\big(\textbf{H}\big(A:\D^2u\big)\big) \ot  \D\big(\textbf{H}\big(A:\D^2u\big)\big)\Big]\, =\, 0
\]
which by decomposing the positive matrix $A$ as $A=A^{1/2}A^{1/2}$, we reformulate as 
\[
\textbf{H}'\big(A:\D^2u\big) \Big|A^{1/2}\D\big(\textbf{H}\big(A:\D^2u\big)\big) \Big|^2\, =\, 0.
\]
By the PDE it follows that on the open set $\Om^*:=\big\{\textbf{H}'\big(A\!:\!\D^2u\big)\neq 0\big\}$ we have $\D\big(\textbf{H}\big(A\!:\!\D^2u\big)\big)=0$ because $A^{1/2}$ is invertible. On the other hand, on $\Om\set\Om^*$ we have $\textbf{H}'\big(A\!:\!\D^2u\big)=0$ which gives $\D\big(\textbf{H}\big(A:\D^2u\big)\big) = \textbf{H}'\big(A\!:\!\D^2u\big)\, \D\big( A\!:\!\D^2u\big)=0$. Consequently, $\D \big(\textbf{H}\big( A\!:\!\D^2u\big)\big)=0$ on the whole of $\Om$ and by the connectivity of the domain we infer that $\textbf{H}\big( A\!:\!\D^2u\big)\equiv c$ for some $c\geq 0$.

\smallskip

\noi (3) $\Rightarrow$ (2): If $\textbf{H}\big(A\!:\!\D^2u\big)\equiv c$ on $\Om$, since  $\{\textbf{H}=c\}$ consists of at most $2$ points and $A\!:\!\D^2u(\Om')$ is a connected set, we obtain $A\!:\!\D^2u(\Om') \sub \{C\}$ where $C\in \{\textbf{H}=c\}$.

\smallskip

\noi (2) $\Rightarrow$ (1): Let us denote the $n$-Lebesgue and the $n-1$-Hausdorff measure as in \cite{EG} by $\mL^n$ and $\mH^{n-1}$ respectively. Fix $\Om'\Subset \Om$ and $\phi\in W^{2,\infty}_0(\Om')$. Extend $\phi$ by zero on $\R^n\set\Om'$, consider the standard mollification $\eta^\e*\D^2\phi$ of it by convolution and let $\Om'_\de$ be a piecewise smooth domain containing the support of $\eta^\e*\D^2\phi$ (e.g.\ union of finite many balls). By the Gauss-Green theorem, we have
\[
\int_{\Om'_\de}\eta^\e*\D^2\phi\,d\mL^n \,=\, \int_{\Om'_\de}\D\big(\eta^\e*\D\phi\big)\,d\mL^n \,=\, \int_{\p\Om'_\de}\big(\eta^\e*\D\phi\big) \ot \nu \,d\mH^{n-1}\, =\, 0
\]
and by letting $\e \ri 0$ and $\de\ri 0$, we get that the average of $\D^2\phi$ over $\Om'$ vanishes. Hence, since $\A\!:\!\D^2u \equiv C$ on $\Om$, we obtain 
\[
C\,=\, \av_{\Om'}A:\D^2u\,d\mL^n\, =\, \av_{\Om'} \Big(A:\D^2u + A:\D^2\phi\Big)\,d\mL^n
\]
By applying $\textbf{H}$ to the above equality, Jensen's inequality for level-convex functions (see e.g.\ \cite{BJW1, BJW2}) implies
\[
\begin{split}
\underset{\Om'}{\ess\,\sup}\,\textbf{H}\big(A:\D^2u\big) \, =\, \textbf{H}(C)
\, &=\, \textbf{H}\left(\av_{\Om'} \Big(A:\D^2u + A:\D^2\phi\Big)\,d\mL^n\right)
\\
&\leq\,\underset{\Om'}{\ess\,\sup}\, \textbf{H}\Big(A:\D^2u + A:\D^2\phi\Big)
\end{split}
\]
which yields that $u\in \mathrm{AM}^2(\E_\infty,\Om)$.

\noi (1) $\Rightarrow$ (4): Proved in part (A) above. The proof of Theorem \ref{theorem5} is now complete. \qed

\ms

\BPCOR \ref{corollary6}. In view of the equivalences among (1)-(4) of Theorem \ref{theorem5}(B), it follows that either the Dirichlet problem for the PDE or for $2$nd order Absolute Minimisers reduces to the uniqueness of solution in $C^3(\Om)\cap W^{2,\infty}_g(\Om)$ to the Dirichlet problem for the linear elliptic $2$nd order PDE
\[
\left\{
\begin{array}{rl}
A: \D^2u \, =\, C,\ \ & \ \text{ in }\Om,
\\
u\,=\,g, \ \D u\,= \, \D g,\ \ & \ \text{ on }\p\Om,
\end{array}
\right.
\]
for some $C\in \R$, which is over-determined and has at most $1$ solution.  \qed
\ms

\section{Existence of $\mD$-solutions to the Dirichlet problem for $\A^2_\infty $} \label{section7}

Herein we establish the existence of $\mD$-solutions with extra properties to the Dirichlet problem for \eqref{1.2}. These solutions are in a sense ``critical points" of \eqref{1.1} and  generally non-minimising and non-unique. They are obtained \emph{without imposing any kind of convexity, neither level-convexity nor quasiconvexity nor ``BJW-convexity" (the notion of $L^\infty$-quasiconvexity of \cite{BJW1})}. Actually, our \emph{only assumption on $\H$ is that it is $C^1$ and depends on $X$ via $X^2=X^\top X$}. 

The method we employ has two main steps. First, given $\Om \sub \R^n$ open and $g \in W^{2,\infty}(\Om)$, we solve the fully nonlinear PDE 
\beq \label{7.1}
\ \ \ \ \ \H\big(\D^2u\big)\, =\, C, \ \ \text{ a.e.\ on }\Om,
\eeq
for admissible large enough ``energy level" $C>0$ depending on the data $g$. For this we use the celebrated \emph{Baire Category method} of Dacorogna-Marcellini (see \cite{DM, D}) which is a convenient analytic alternative to Gromov's Convex Integration. Next, we use the machinery of $\mD$-solutions to make the next \emph{non-rigorous statement} precise: every solution $u$ to \eqref{7.1} solves \eqref{1.2} because $\D\big(\H\big(\D^2u\big)\big)\equiv 0$ and \eqref{1.2} ``equals" \eqref{1.2a}. This is indeed true in the class of classical/strong solutions in $C^3(\Om)$ or $W^{3,\infty}(\Om)$, but not in the natural $W^{2,\infty}(\Om)$ class. This method of constructing critical point solutions has previously been applied successfully to the vector-valued first order case and its generalisations, see \cite{K8, K9, AK, CKP}. The principal result of this section therefore is:

\begin{theorem}[Existence of $\mD$-solutions to the Dirichlet problem for $\A^2_\infty u =0 $] \label{theorem12a} 

Let $\H :\R^{n^{\ot 2}}_s\larrow \R$ be such that $\H(X) = h\big(X^2\big)$ for some $h \in C^1\big(\R^{n^{\ot 2}}_s\big)$. Consider also an open set $\Om\sub \R^n$ and fix $g\in W^{2,\infty}(\Om)$. Then, for any ``energy level" $c>\|\D^2g\|_{L^\infty(\Om)}$, the Dirichlet problem for \eqref{1.2}
\beq \label{7.2}
\left\{
\begin{array}{rl}
\A^2_\infty u \, =\, 0, & \ \text{ in }\Om, \smallskip\\
u\, =\, g, \ \D u \, =\, \D g, & \ \text{ on }\p\Om,
 \end{array}
 \right.
\eeq
has (an infinite set of) $\mD$-solutions in the class
\[
\mA_c\,:=\, \Big\{v\in W^{2,\infty}_g(\Om)\ :\ \H\big(\D^2v\big) = \H\big(c\mathrm{I}\,\big) \, \text{ a.e.\ on }\Om \Big\}.
\]
Namely, there is a set of u's in $\mA_c$ such that (in view of Definitions \ref{definition2} and  \ref{definition3})
\beq \label{7.3}
\ \ \ \ \ \ \ \int_{\smash{\overline{\R}}^{n^{\ot 3}}_s} \Phi(\X)\,\left[\H_X\big(\D^2u(x)\big)^{\ot 3}\!: \X^{\,\ot 2}\right]\, d\big[\mD^3 u(x) \big](\X)\, =\, 0, \ \ \text{ a.e.\ }x\in \Om,
\eeq
for any diffuse $3$rd derivative $\mD^3 u \in \mY\big(\Om,\smash{\overline{\R}}^{n^{\ot 3}}_s\big)$ and any $\Phi \in C_c\big( {\R}^{n^{\ot 3}}_s \big)$.

\smallskip
\end{theorem}

\BPT \ref{theorem12a}. Let $\Om\sub \R^n$ be a given open set and $g\in W^{2,\infty}(\Om)$, $n\in \N$. We begin by showing the next result.

\begin{claim} \label{claim13} For any fixed $c>\|\D^2 g\|_{L^\infty(\Om)}$, there exist (an infinite set of) solutions in $ W_g^{2,\infty}(\Om)$ such that
\[
\left\{\ \ \
\begin{array}{rl}
\D^2u^\top \D^2 u \, =\, c^2\mathrm{I},  & \ \text{ a.e.\ in }\Om, 
\smallskip\\
u\,=\,g, \ \D u\,=\,\D g, & \ \text{ on }\p\Om.
 \end{array}
 \right.
\]
\end{claim}

\BPC \ref{claim13}. Let $\{\la_1(X),...,\la_n(X)\}$ symbolise the eigenvalues of the symmetric matrix $X\in \R^{n^{\ot 2}}_s$ in increasing order. By the results of [\cite{DM}, p.\ 200] the Dirichlet problem
\[
\left\{\ \
\begin{array}{rl}
\big|\la_i\big(\D^2 v\big)\big|\, =\, 1, \ \ \ & \text{ a.e.\ on }\Om, \ \al=1,...,n, \smallskip
\\
v\,=\,g/c, & \text{ on }\p\Om,
 \end{array}
 \right.
\]
has solutions $v \in W_{g/c}^{2,\infty}(\Om,\R^N)$ because we have that
\[
\begin{split}
\max_{i=1,...,n} \left\{\underset{\Om}{\ess\,\sup}\, \Big| \la_i\Big(\D^2\big(g/c\big)\Big) \Big|
\right\} &\leq \, \frac{1}{c}\, \underset{\Om}{\ess\,\sup}\left\{ \max_{|e|=1}\Big|  \D^2 g : (e\ot e) \Big|\right\}
\\ 
&\leq\, \frac{1}{c}\, \big\| \D^2 g\big\|_{L^\infty(\Om)} \, <\, 1,
\end{split}
\]
a.e.\ on $\Om$. By rescaling as $u:=cv$, we get existence of solutions $u \in W_{g}^{2,\infty}(\Om)$ to
\[
\left\{ \ \
\begin{array}{rl}
\big|\la_i\big(\D^2 u\big)\big|\, =\, c, \, & \text{ a.e.\ on }\Om, \ \al=1,...,n, \smallskip
\\
v\,=\,g, & \text{ on }\p\Om.
 \end{array}
 \right.
\]
Note now that $\la_i\big(\D^2 u\big)^{\!2}=\la_i\big(\D^2 u^\top\D^2u\big)$ for all $i=1,...,n$. Hence, by the Spectral Theorem there is an $L^\infty$ map with values in the orthogonal matrices $O : \Om\sub \R^n \larrow O(n,\R) \sub \R^{n^{\ot 2}}$ such that
\[
\D^2u^\top  \D^2u\, =\, O 
\left[
\begin{array}{lr}
\la_1\big(\D^2 u\big)^2 &   \textbf{0}\
\\
\hspace{60pt}  \ddots &
\\
\ \textbf{0}    & \la_n\big(\D^2 u\big)^2
\end{array}
\right]
O^\top\, =\, O\, (c^2\mathrm{I})\, O^\top\, =\, c^2\mathrm{I},
\]
a.e.\ on $\Om$. The claim thus ensues.   \qed
\ms

Now we complete the proof of the theorem. By our assumption on $\H$, for any $u \in W_{g}^{2,\infty}(\Om)$ as in Claim \ref{claim13} we have
\beq  \label{7.4}
\H \big(\D^2 u\big)\, =\, h\big(\D^2u^\top \D^2 u\big)\, =\, h\big(c^2\mathrm{I}\big)\, =\, h\big((c\mathrm{I})^2\big)\, =\,\H \big( c\mathrm{I}\big),
\eeq
a.e.\ on $\Om$. Hence, $u\in \mA_c$. Note also that by \eqref{7.4} we have $\H\big(\D^2 u(x)\big)=\, $const for a.e.\ $x\in \Om$. The next claim completes the proof.

\begin{claim} 
\label{claim} 
If $\H\big(\D^2 u \big)=C$ a.e.\ on $\Om$, then $u$ is a $\mD$-solution to \eqref{1.2}, that is \eqref{7.3} holds for any fixed $\Phi \in C_c\big( {\R}^{n^{\ot 3}}_s \big)$ and a.e.\ $x\in \Om$ for any diffuse $3$rd derivative of $u$
\beq \label{7.5}
\ \ \ \ \ \ \de_{\D^{1,h_{m}}\D^2 u}\weakstar \,\mD^3 u \  \ \text{ in }\mY\big(\Om,\smash{\overline{\R}}^{n^{\ot 3}}_s\big), \ \ \text{ as $m\ri \infty$}.
\eeq 
\end{claim}

\BPC \ref{claim}. Fix such an $x\in \Om$, $0<|h|<\dist(x,\p\Om)$ and $k\in \{1,...,n\}$. By Taylor's theorem, we have
\[
\begin{split} 
 0\, &=\, \H\big(\D^2 u(x+he^k)\big)-\H\big(\D^2 u(x)\big)
\\
&=\,\sum_{i,j}\int_0^1 \H_{X_{i j}}\Big(\D^2 u(x) + \la \Big[\D^2 u(x+he^k)-\D^2 u(x)\Big] \Big)\,d\la\, \centerdot
\\ 
& \hspace{40pt}  \centerdot \Big[\D^2_{ij} u(x+he^k)-\D^2_{ij} u(x)\Big].
\end{split}
\]
This implies for any $k=1,...,n$ the identity
\[
\begin{split} 
    \sum_{i,j} & \, \H_{X_{i j}}\big(\D^2 u(x)\big)\, \big(D^{1,h}_k\D^2_{ij}u \big)(x) 
\, =
\end{split}
\]
\beq \label{7.6}
\begin{split}
 &=\,  \sum_{i,j} \bigg\{\int_0^1 \bigg[-\, \H_{X_{i j}}\Big(\D^2 u(x)  + \la \Big[\D^2 u(x+he^k)-\D^2 u(x)\Big] \Big)
\\
& \hspace{62pt} +\, \H_{X_{i j}}\big(\D^2 u(x)\big)\bigg]\,d\la \bigg\} \big(D^{1,h}_k\D^2_{ij}u \big)(x)
\\
&=:\ \sum_{ij}\mathcal{E}_{ijk}(x,h)\, \big(D^{1,h}_k\D^2_{ij}u \big)(x)
\end{split}
\eeq
where $\mathcal{E}_{ijk}$ is the ``error tensor". By taking \eqref{7.6} for $k=p,q$, multiplying these two equations with $\H_{X_{p q}}\big(\D^2 u(x)\big)$ and summing in $p,q\in \{1,...,n\}$, we obtain
\beq \label{7.7}
\begin{split} 
    \sum_{i,j,r,s,p,q}& \Bigg[ \bigg( \H_{X_{p q}}\big(\D^2 u\big) \, \H_{X_{i j}}\big(\D^2 u\big)  \, \H_{X_{rs}}\big(\D^2 u\big)\bigg)  \big(D^{1,h}_p\D^2_{ij}u \big)  \big(D^{1,h}_q\D^2_{rs}u \big)\Bigg]
\\
& =\, \sum_{i,j,r,s,p,q} \bigg\{ \H_{X_{p q}}\big(\D^2 u\big)  \, \mathcal{E}_{ijp}(\cdot,h) \, \mathcal{E}_{rsq}(\cdot,h)\bigg\} \big(D^{1,h}_p\D^2_{ij}u \big)  \big(D^{1,h}_q\D^2_{rs}u \big) 
\\
&=:\, \sum_{i,j,r,s,p,q}\ \mathcal{E}_{pqijrs}(\cdot ,h)\, \big(D^{1,h}_p\D^2_{ij}u \big)  \big(D^{1,h}_q\D^2_{rs}u \big) .
\end{split}
\eeq
Let  $(h_m)_1^\infty$ be an infinitesimal sequence giving rise to a diffuse $3$rd derivative as in \eqref{7.5}. We rewrite \eqref{7.7} for $h=h_m$ compactly as 
\[
\ \ \bigg(\big(\H_X\big(\D^2u\big)\big)^{\ot 3} \, -\  \mathcal{E}(\cdot,h_m)\bigg):\big(\D^{1,h_m}\D^2u\big)^{\ot 2}\,=\,0 
\]
for $m\in \N$. Then for any $\Phi \in C_c\big( {\R}^{n^{\ot 3}}_s \big)$ and $\phi \in C_c(\Om)$, this yields
\beq \label{7.8}
\ \ \int_\Om\phi \int_{ \smash{\overline{\R}}^{n^{\ot 3}}_s } \Phi(\X)\, \bigg[\Big(\big(\H_X\big(\D^2u\big)\big)^{\ot 3}  -\  \mathcal{E}(\cdot,h_m) \Big): \X^{\ot 2}\bigg] \, d\big[ \de_{\D^{1,h_m}\D^2u} \big](\X)\,=\,0 .
\eeq
Since $|\D^2u| \in L^{\infty}(\Om)$, by the continuity of the translation operation in $L^1$ we have $|\D^2u(\cdot + z) -\D^2u| \larrow 0$ as $z\ri 0$, in $L^1_{\text{loc}}(\Om)$. Hence, along perhaps a further subsequence $(m_i)_{i=1}^\infty$ we have $\D^2 u\big(x+h_m e^k\big) \larrow \D^2u(x)$ for a.e.\ $x\in \Om$ as $m\ri \infty$, $k=1,...,n$. Since $\H \in C^1\big(\R^{n^{\ot 2}}_s\big)$ and $|\D^2 u| \in L^{\infty}(\Om)$, the Dominated Convergence theorem and the definition of the errors in \eqref{7.6}-\eqref{7.7} imply that $\big|\mathcal{E}(\cdot,h_m)\big| \larrow 0$ in $L^1_{\text{loc}}(\Om)$ subsequentially as $m\ri \infty$. We define the Carath\'eodory functions
\[
\ \ \ 
\begin{split}
\Psi^m(x,\X)\, &:=\, \phi(x) \, \Phi(\X)\Big| \Big(\big(\H_X\big(\D^2u\big)\big)^{\ot 3}  -\  \mathcal{E}(\cdot,h_m) \Big): \X^{\ot 2}  \Big|,
\\
\Psi^\infty(x,\X)\, &:=\, \phi(x)\, \Phi(\X) \Big| \big(\H_X\big(\D^2u\big)\big)^{\ot 3} : \X^{\ot 2} \Big|
\end{split}
\]
which are elements of the Banach space $L^1\big(E,C\big( \smash{\overline{\R}}^{n^{\ot 3}}_s  \big)\big)$ (because of the compactness of the supports of $\phi,\Phi$) and we also have
\beq \label{7.10}
\Psi^{m} \larrow \Psi, \ \text{ as }m \ri \infty \ \text{ in }L^1\big(E,C\big( \smash{\overline{\R}}^{n^{\ot 3}}_s \big)\big)
\eeq
which is a consequence of that $\big|\mathcal{E}(\cdot,h_m)\big| \larrow 0$ in$L^1_{\text{loc}}(\Om)$ and the estimate
\[
\|\Psi^{m_k} -\Psi^\infty\|_{L^1(E,C( \smash{\overline{\R}}^{n^{\ot 3}}_s ))}\, \leq\, \max_{\X \in \supp(\Phi)}\big|\Phi(\X)\big||\X|^2 \int_{\supp(\phi)} \big|\phi \,\mathcal{E}(\cdot,h_m)\big| .
\]
The weak*-strong continuity of the duality pairing \eqref{dp}, \eqref{7.10} and \eqref{7.5} allow us to pass to the limit in \eqref{7.8} as $m\ri \infty$  and deduce
\[
\ \ \ \int_\Om\phi \int_{ \smash{\overline{\R}}^{n^{\ot 3}}_s } \Phi(\X)\, \bigg[ \big(\H_X\big(\D^2u\big)\big)^{\ot 3} : \X^{\ot 2}\bigg] \, d [\mD^3u ](\X)\,=\,0 .
\]
Since $\phi\in C_c(\Om)$ is arbitrary, \eqref{7.3} follows and the claim ensues. \qed
\ms

The proof of the theorem is now complete. \qed

\section{Explicit $p$-Biharmonic and $\infty$-Biharmonic functions in $1D$}  \label{section8}

In this section we give explicit solutions to the Dirichlet problem for the $p$-Bilaplacian and the $\infty$-Bilaplacian when $n=1$. In this case the equations are
\begin{align}
 \De^2_p u\, =\, \big(|u''|^{p-2}u''\big)''\,&=\,0  \label{8.1}
 \\
 \De^2_\infty u\,=\, (u'')^3(u''')^2\, &=\, 0   \label{8.2}
\end{align}
The weak solutions we construct for \eqref{8.1} are obtained by solving the equation explicitly  for even exponents $p\in 2\N$, while the $\mD$-solutions we construct for \eqref{8.2} are piecewise quadratic. In either case the solutions in general have at least $1$ singular point in their domain, unless the boundary data can be interpolated by a quadratic polynomial function in which case the solutions are smooth. Accordingly, the main result of this section is: 

\begin{theorem}[Explicit generalised solutions] \label{theorem12} Consider $a,b,A,B,A',B' \in \R$ with $a<b$ and set
\[
E\,:=\, \frac{B'-A'}{b-a}\,-\,\frac{2(B -A -A'(b-a))}{(b-a)^2}.
\]

\noi \emph{(A) [$p$-Bilaplacian]} Let $p\in 2\N$. Then, the problem 
\beq  \label{8.3}
\left\{\ \ 
\begin{split}
&\De^2_p u\, =\,0, \text{ in }(a,b)\sub \R,
\ms\\
& u(a)\, =\,A,\ \ u(b)\, =\,B, 
\ms \\
& u'(a)\, =\,A',\ \ u'(b)\, =\,B', 
\end{split}
\right.
\eeq
has a unique weak solution $u_p\in W^{2,p}(a,b)$ given by:

\noi  \emph{(i)} For critical data satisfying $E=0$ (which can be interpolated by a quadratic function),  
\beq  \label{8.4}
u_p(x)\,=\,A+A'(x-a)+ \frac{1}{2}\left( \frac{B'-A'}{b-a} \right)(x-a)^2.
\eeq

\noi \emph{(ii)} For non-critical data satisfying $E\neq0$ (which can not be interpolated by a quadratic function),  
\beq  \label{8.5}
u_p(x)\,=\,A+ \left(A' -\frac{p-1}{p\la} |\la a +\mu|^{\frac{p}{p-1}} \right)(x-a)+ 
\frac{p-1}{p\la} \int_a^x|\la t +\mu|^{\frac{p}{p-1}}\,dt.
\eeq
Further, $(\la,\mu) \in (\R\set\{0\})\by \R$ is the unique solution to the algebraic equations:
\beq  \label{8.6}
\left\{
\begin{split}
 (B'-A')\frac{p\la}{p-1}\, &=\,  |\la b +\mu|^{\frac{p}{p-1}}-|\la a +\mu|^{\frac{p}{p-1}}, 
\\
 \big(B-A-A'(b-a)\big)\frac{p\la}{p-1}\, &=\,  \int_a^b|\la t +\mu|^{\frac{p}{p-1}}dt-|\la a +\mu|^{\frac{p}{p-1}} (b-a).
\end{split}
\right.
\eeq
In particular, $u_p \in C^\infty\big((a,b)\set \{-\mu/\la\}\big)$.

\ms

\noi \emph{(B) [$\infty$-Bilaplacian]} For any large enough ``energy level" $C>0$ depending only on the boundary data (see \eqref{8.11}), the problem \eqref{8.3} for $p=\infty$ has a unique piecewise quadratic $\mD$-solution $u_\infty\in W^{2,\infty}(a,b)$ given by 
\beq  \label{8.8}
u_\infty(x)\, =\, A+A'(x-a)+ C\int_a^x\Big[\mL^1\big([a,t]\cap I_C\big)-\mL^1\big([a,t]\set I_C\big)\Big]\,dt.
\eeq
Here $I_C=[x_C,y_C]$ is the interval with endpoints
\beq  \label{8.9}
x_C\,=\, \frac{-K-L^2+2bL}{2L}, \ \ \ y_C\, =\, \frac{-K+L^2+2bL}{2L}
\eeq
where
\beq  \label{8.10}
K\,=\, \frac{2\big(B-A-A'(b-a)\big) + C(b-a)^2}{2C},\ \ \ L\,=\, \frac{B'-A' +C(b-a)}{2C}.
\eeq
In particular, $u_\infty$ satisfies $|u_\infty''|=C$ a.e.\ on $\Om$ and $u_\infty \in C^\infty\big((a,b)\set\{x_C,y_C\}\big)$.
\end{theorem}

\begin{remark} We note that the solution $u_\infty$ above is not the limit of $u_p$ as $p\ri\infty$. The function $\lim_p u_p$ is indeed Absolutely Minimising by the results of Section \ref{section4} but \emph{we do not prove here that is solves in the $\mD$-sense the equation}. Instead, we solve \eqref{8.2} by solving the fully nonlinear equation $|u''|=C$ for a fixed energy level $C>0$ and using the previous section to characterise it as a $\mD$-solution to \eqref{8.2}. The numerics of the next section show that  $\lim_p u_p$ has at most $1$ ``breaking point" for the $2$nd derivative in the domain of definition, while these solutions are ``critical points" and as such have instead less regularity and $2$ ``breaking points" of their $2$nd derivative.
\end{remark}

\BPT \ref{theorem12}. (A) Let $u_p$ be a weak solution to \eqref{8.3}. By standard convexity and variational arguments (see e.g.\ \cite{E, D}), the solution exists and it is energy minimising and unique. Note now that the function $\R \ni t\mapsto |t|^{p-2}t=t^{p-1} \in \R$ and its inverse $t\mapsto t^{1/(p-1)}$ are odd because $p\in 2\N$. We obtain (i)-(ii) directly by differentiating twice the explicit formulas \eqref{8.4}-\eqref{8.5}. By the previous observation, in either case this gives $u_p''(x) =\big( \la x\, +\, \mu\big)^{\frac{1}{p-1}}$ for $a<x<b$, where in the case of (i) we have $\mu=0$ and $\la=((B'-A')/(b-a))^{p-1}$, whilst in the case of (ii) the parameters $(\la,\mu)$ are given \eqref{8.6}. The latter is just a compatibility condition arising by the boundary conditions. In both cases we get that the function $|u_p''|^{p-2}u_p''$ is affine and $\big(|u_p''|^{p-2}u_p''\big)(x) =\la x +\mu$ for $a<x<b$. As a consequence, $\De^2_p u_p=0$ weakly on $(a,b)$.

\ms

\noi (B) Let $x_C,y_C$ be given by \eqref{8.9}-\eqref{8.10} for any $C>C^*$, $C^*$ large enough to be specified. The compatibility requirement $a\leq x_C < y_C \leq b$ (in order to have endpoints of an interval $I_C=[x_C,y_C]$ which lies inside $[a,b]$) after a calculation is equivalent to
\beq \label{8.11}
\left\{
\begin{split}
\frac{1}{C}\left[ \frac{(B'-A')^2}{4C}+ \frac{(B'-A')(b-a)}{2} -\big(B-A-A'(b-a)\big)\right] \leq \frac{(b-a)^2}{4},
\\
\frac{1}{C}\left[ \frac{(B'-A')^2}{4C} -\frac{(B'-A')(b-a)}{2}\right] \leq \frac{(b-a)^2}{4},
\end{split}
\right.
\eeq
and hence for $C^*$ large enough we indeed have existence. Thus, $u_\infty$ as given by \eqref{8.8} is a well-defined $W^{2,\infty}(a,b)$ function. By differentiating \eqref{8.8}, we get $u_\infty'' = C$ on $(x_C,y_C)$ and $u_\infty'' = -C$ on $(a,x_C)\cup (y_C,b)$, giving $|u_\infty''|=C$ on $(a,b)\set\{x_C,y_C\}$. By Claim \ref{claim} of Section \ref{section7}, $u_\infty$ is a $\mD$-solution to $\De^2_\infty u=0$.  Note now that $x_C,y_C$ satisfy the identities $(x_C)^2-(y_C)^2+2b(y_C-x_C)=K$
and $y_C-x_C = L$ and by a calculation, it can be verified that the algebraic equations above are equivalent to
\[
\left\{
\begin{split}
\ \ \ B \, &=\, A \, +\, A'(b-a)\, + C\int_a^b \Big[\mL^1\big([a,t]\cap I_C\big)-\mL^1\big([a,t]\set I_C\big)\Big]\,dt,
\\
B'\, &=\, A' \,+\, C\,\Big[2\,\mL^1(I_C)-(b-a)\Big].
\end{split}
\right.
\]
The latter pair of equations are just a restatement of the fact that $u_\infty$ satisfies the boundary conditions.  The theorem ensues. \qed
\ms

\section{Numerical approximations of $\infty$-Biharmonic functions}  \label{section9}

In this section we illustrate some of the properties of $\infty$-Biharmonic functions using numerical techniques. We present results from a numerical scheme that makes use of a $p$-Biharmonic approximation, that is, we make use of the derivation through the $p$-limiting process given in Section \ref{section3}. Our numerical scheme of choice is a finite element method and is fully described in \cite{KP2} where we prove for fixed $p$ that the scheme converges to the weak solution of the $p$-Biharmonic problem. The results there illustrate that for practical purposes, as one would expect, the approximation of $p$-Biharmonic functions for large $p$ gives good resolution of candidate $\infty$-Biharmonic functions. In this work for brevity we restrict ourselves to presenting only some results.

\ms
\ms

\noi \textbf{Test 1: the $1$-dimensional problem.} We consider the Dirichlet problem \eqref{8.3} for the $p$-Bilaplacian \eqref{8.1} for $n=1$ with the data $A,B,A',B'$ being given by the values of the cubic function 
\begin{equation} \label{9.2}
  g(x) = \tfrac{1}{120} (4x-3)(2x-1)(4x-1)
\end{equation}
on $(0,1)$. We simulate the $p$-Bilaplacian \eqref{8.1} for increasing values of $p$ and present the results in Figure \ref{fig:floppydonkeydick} indicating that in the limit the $\infty$-Biharmonic function should be piecewise quadratic.

\providecommand{\figwidth}{\textwidth}
\newcommand{\figscale}{1.1}

{

\begin{figure}[!ht]
  \caption[]
  {\label{fig:floppydonkeydick}
    A mixed finite element approximations to an $\infty$-Biharmonic function using $p$-Biharmonic functions for various $p$ for the problem given by \eqref{8.3} and \eqref{9.2}. Notice that as $p$ increases, $u''$ tends to a piecewise constant up to Gibbs oscillations. This is an indication the solution is indeed piecewise quadratic. 
  }
  \begin{center}
   \subfigure[{\label{fig:a1}
        The approximation to $u''$, the Laplacian of the solution of the $4$-Bilaplacian.
    }]{
      \includegraphics[scale=\figscale,width=0.46\figwidth]{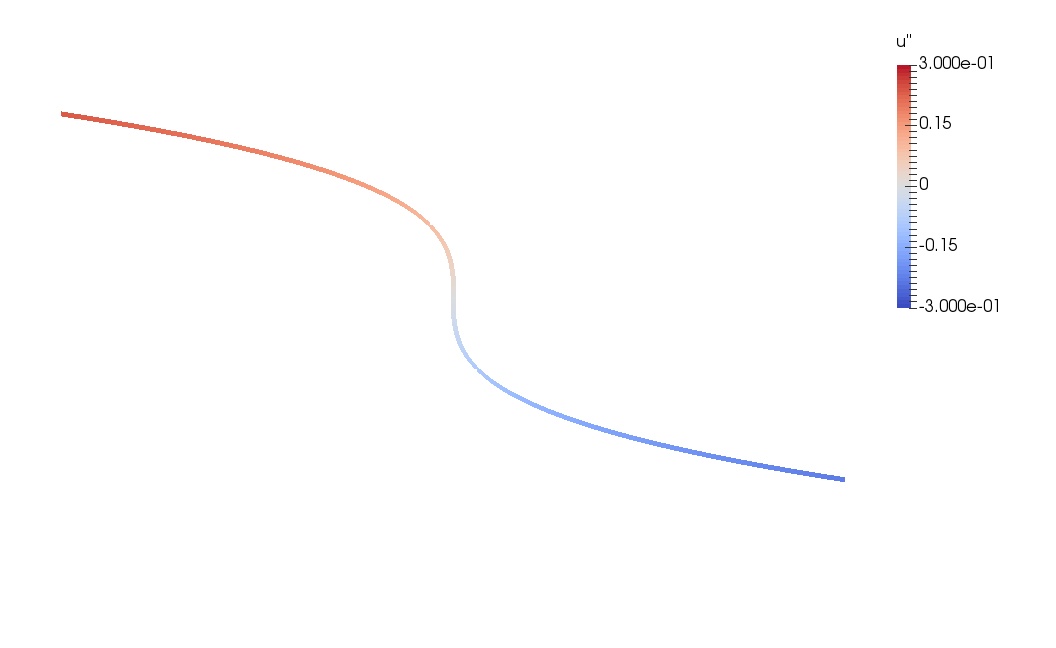}
    }
    \hfill
    \subfigure[{\label{fig:a2}
                The approximation to $u''$, the Laplacian of the solution of the $12$-Bilaplacian.
    }]{
      \includegraphics[scale=\figscale,width=0.46\figwidth]{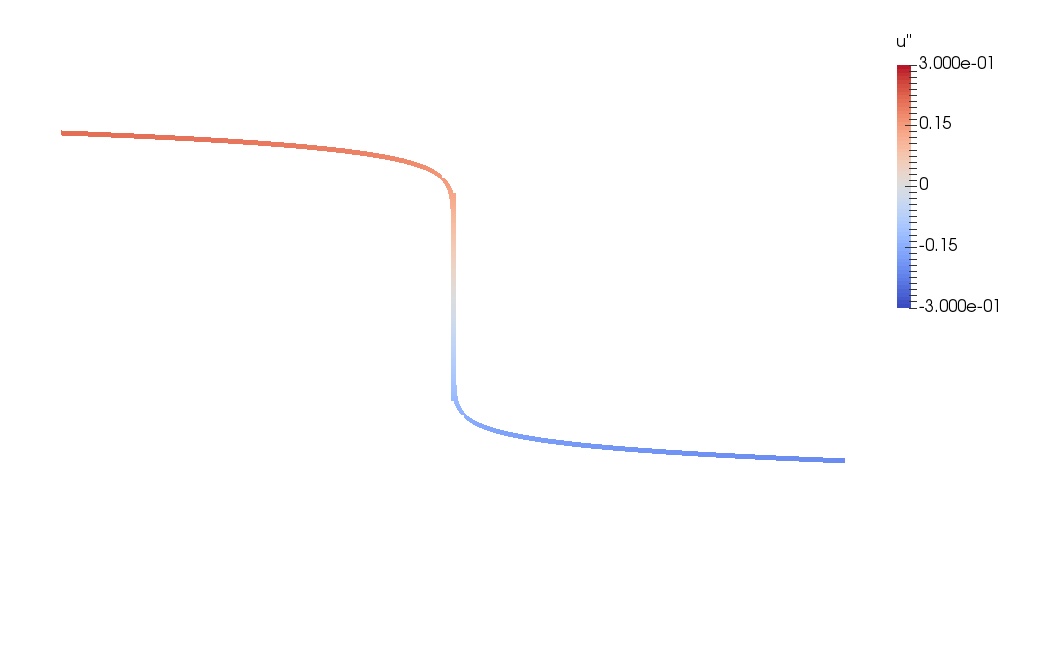}
    }
    \hfill
    \subfigure[{\label{fig:a2}
                The approximation to $u''$, the Laplacian of the solution of the $42$-Bilaplacian.        
    }]{
      \includegraphics[scale=\figscale,width=0.46\figwidth]{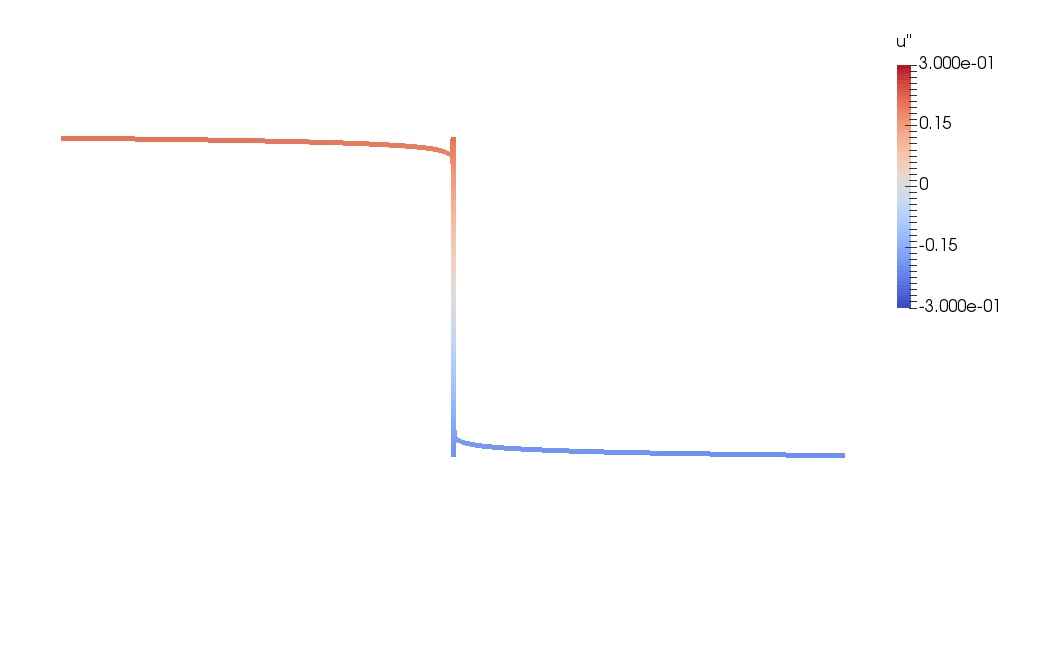}
    }
    \hfill
    \subfigure[{\label{fig:a2}
                        The approximation to $u''$, the Laplacian of the solution of the $202$-Bilaplacian.
    }]{
      \includegraphics[scale=\figscale,width=0.46\figwidth]{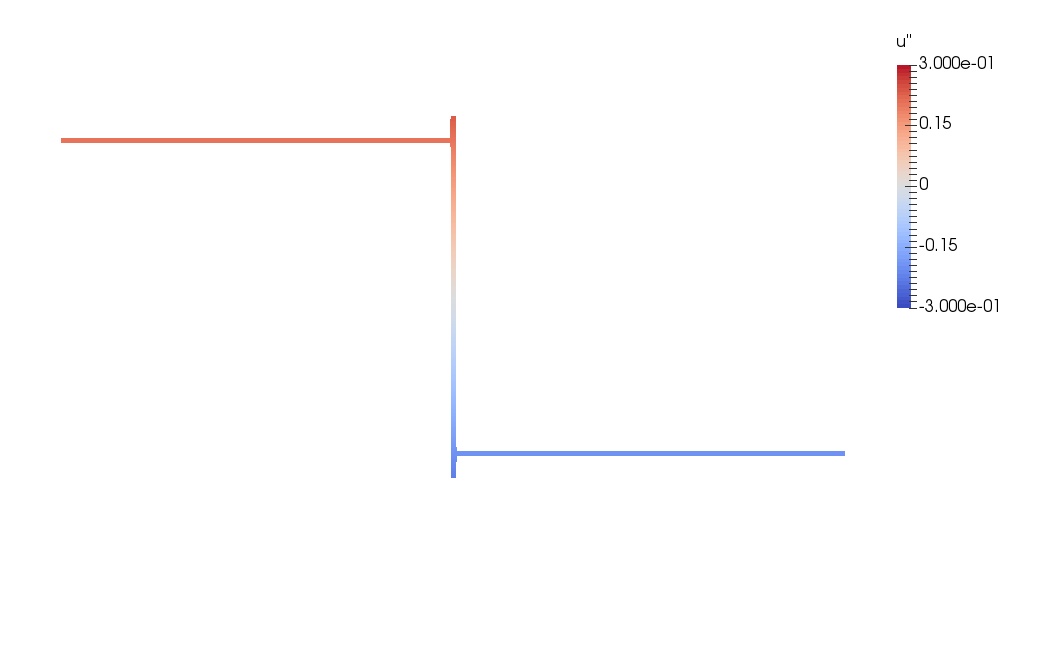}
    }
    \subfigure[{\label{fig:a1}
        The approximation to $u$, the solution of the $4$-Bilaplacian.
    }]{
      \includegraphics[scale=\figscale,width=0.46\figwidth]{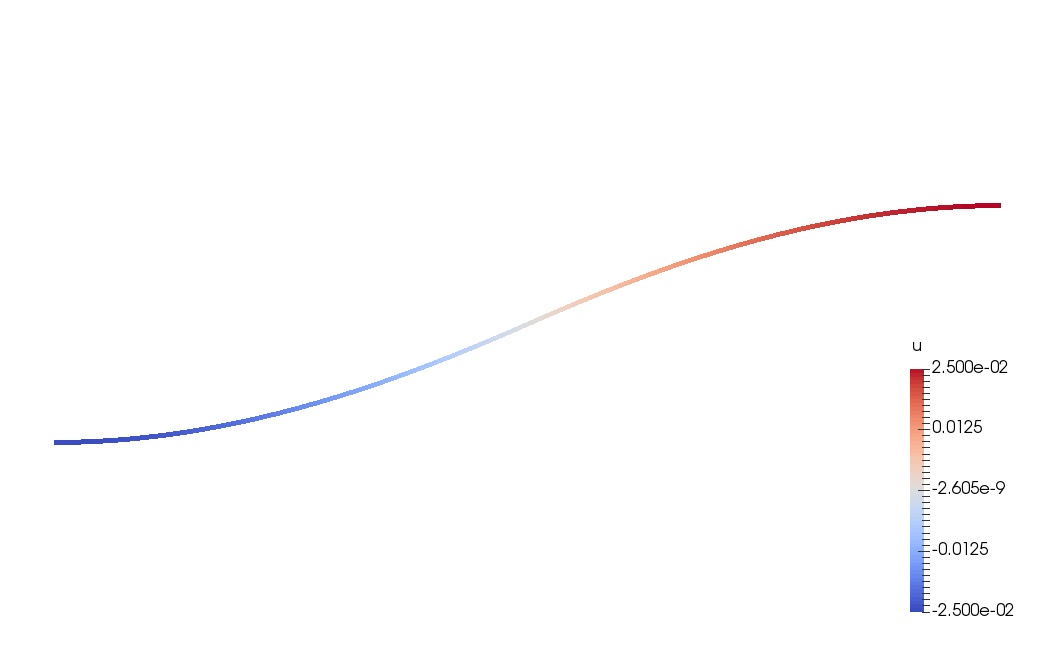}
    }
    \hfill
    \hfill
    \subfigure[{\label{fig:a2}
                The approximation to $u$, the solution of the $202$-Bilaplacian.
    }]{
      \includegraphics[scale=\figscale,width=0.4\figwidth]{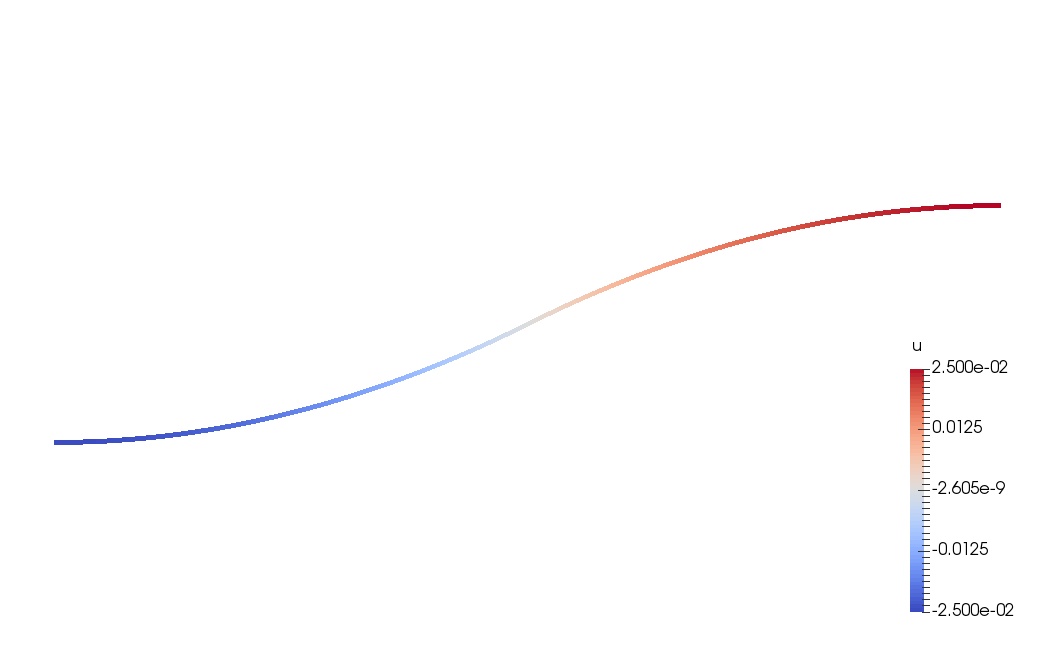}
    }    
  \end{center}
  \end{figure}

}

\ms
\ms

\noi \textbf{Test 2: the $2$-dimensional problem.} Now we illustrate some of the complicated behaviour of the $p$-Bilaplacian for $n=2$:
\begin{equation}
  \label{eq:2dpbiharm}
  \left\{ \ \ \
  \begin{array}{rl}	
    \Delta \big( |\Delta u|^{p-2} \Delta u \big) \,=\,0, \ \ & \text{ in }\Om = [-1,1]^2,
    \smallskip\\
    u \,= \,g, \ \D u \,=\, \D g, \ \ &   \text{ on }\p\Om,
  \end{array}
  \right.
\end{equation}
where $g$ is prescribed as
\begin{equation} \label{9.4}
  g(x,y) = \tfrac{1}{20} \cos(\pi x) \cos(\pi y).
\end{equation}
We simulate the $p$-Bilaplacian for increasing values of $p$ and present the results in Figure \ref{fig:floppydonkeydick3} indicating that in the limit the $\infty$-Biharmonic function should be piecewise quadratic however the behaviour is quite unexpected and complicated interface patterns emerge even with this relatively simple boundary data.

{

\begin{figure}[!ht]
  \caption[]
  {\label{fig:floppydonkeydick3}
    A mixed finite element approximations to an $\infty$-Biharmonic function using $p$-Biharmonic functions for various $p$ for the problem given by \eqref{eq:2dpbiharm} and \eqref{9.4}. Notice that as $p$ increases, $\Delta u$ tends to be piecewise constant. This is an indication the solution satisfies the Poisson equation with piecewise constant right hand side albeit with an extremely complicated solution pattern that clearly warrants further investigation.
  }
  \begin{center}
    \subfigure[{\label{fig:a1}
        The approximation to $\Delta u$, the Laplacian of the solution of the $4$-Bilaplacian.
    }]{
      \includegraphics[scale=\figscale,width=0.47\figwidth]{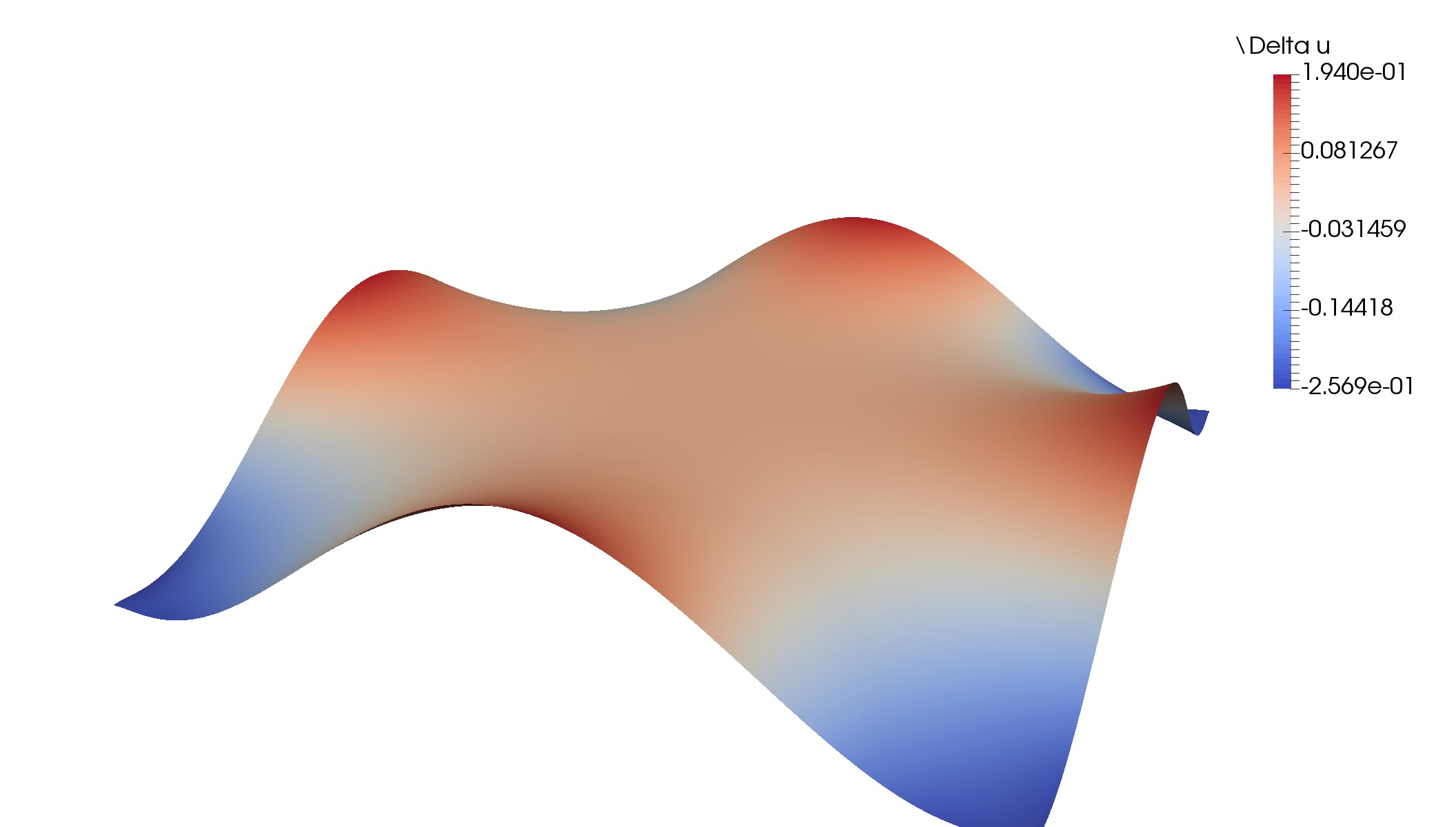}
    }
    \hfill
    \subfigure[{\label{fig:a2}
        The approximation to $\Delta u$, the Laplacian of the solution of the $42$-Bilaplacian.        
    }]{
      \includegraphics[scale=\figscale,width=0.47\figwidth]{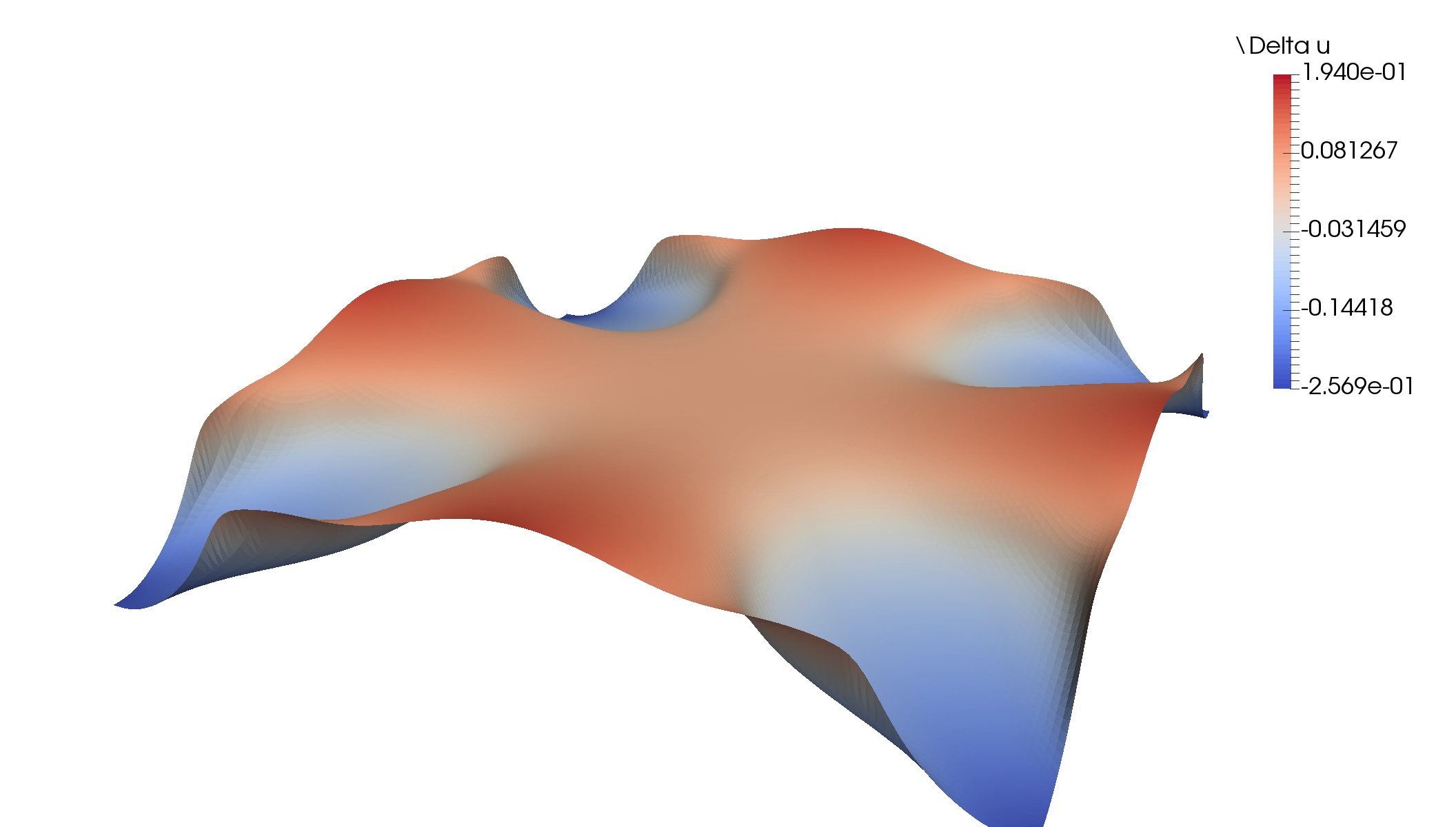}
    }
    \hfill
    \subfigure[{\label{fig:a2}
        The approximation to $\Delta u$, the Laplacian of the solution of the $68$-Bilaplacian.
    }]{
      \includegraphics[scale=\figscale,width=0.47\figwidth]{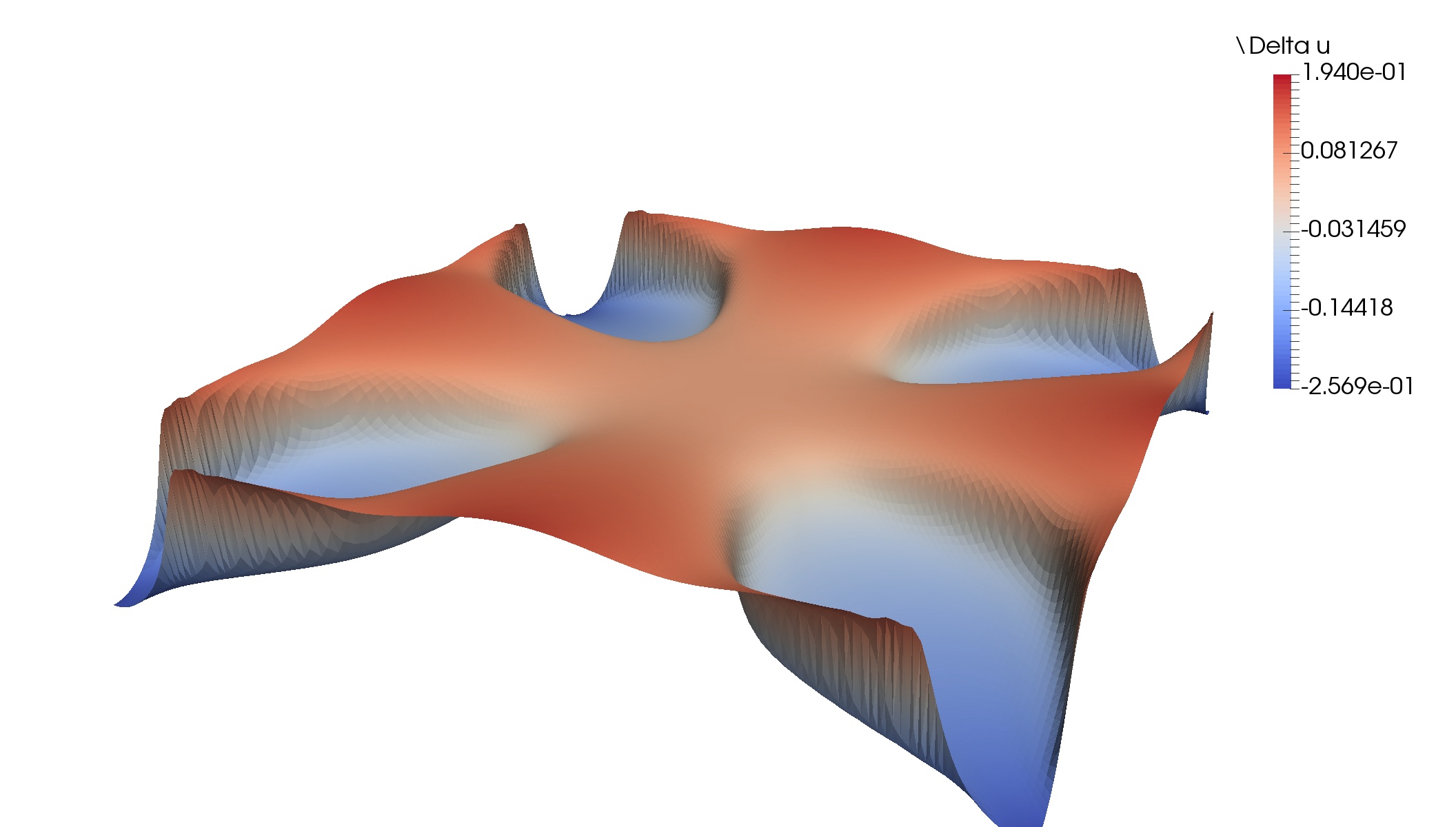}
    }
    \hfill
    \subfigure[{\label{fig:a2}
        The approximation to $\Delta u$, the Laplacian of the solution of the $142$-Bilaplacian.
    }]{
      \includegraphics[scale=\figscale,width=0.47\figwidth]{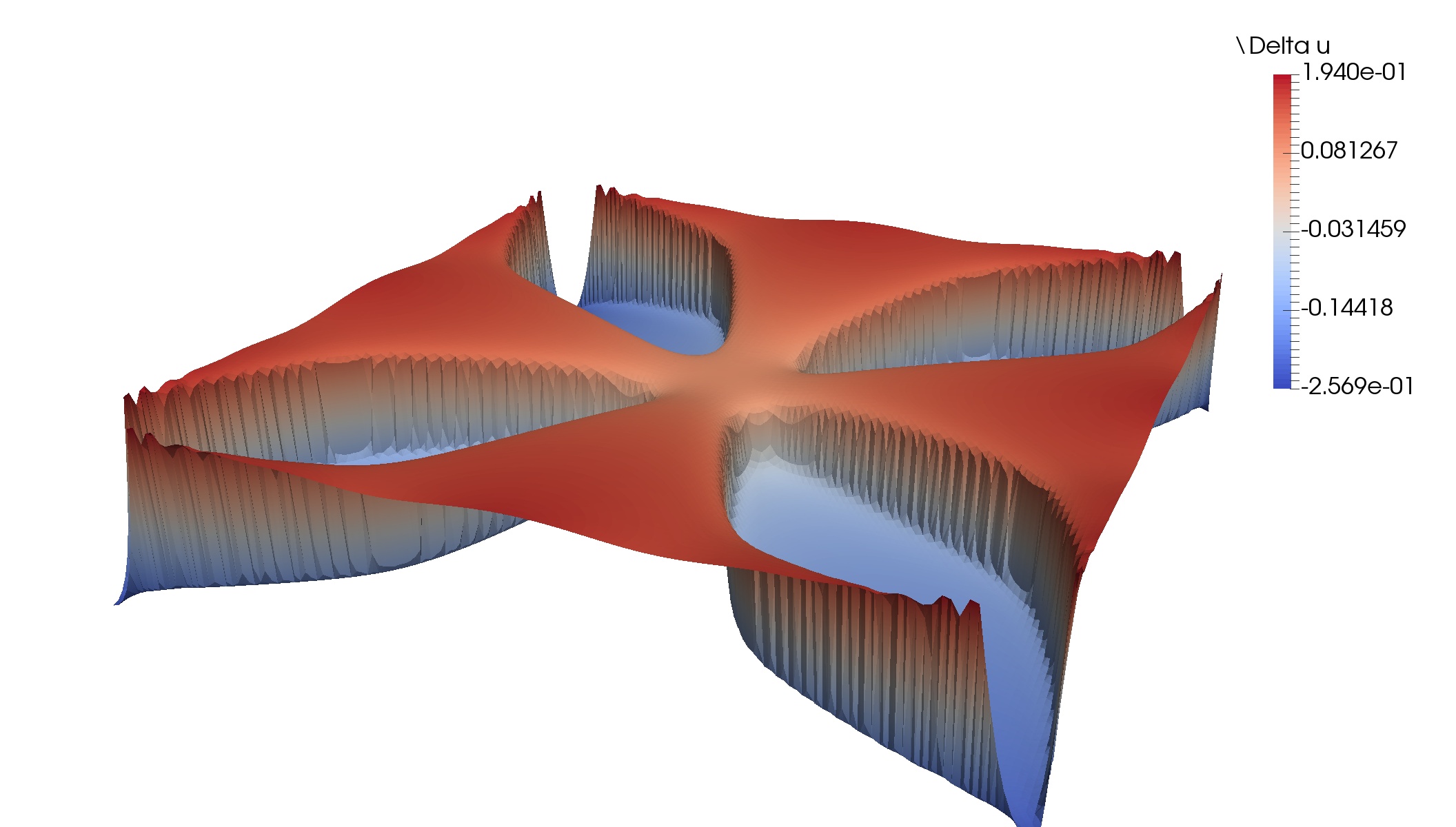}
    }
 
    \hfill
    \subfigure[{\label{fig:a1}
        The approximation to $u$, the solution of the $4$-Bilaplacian.
    }]{
      \includegraphics[scale=\figscale,width=0.47\figwidth]{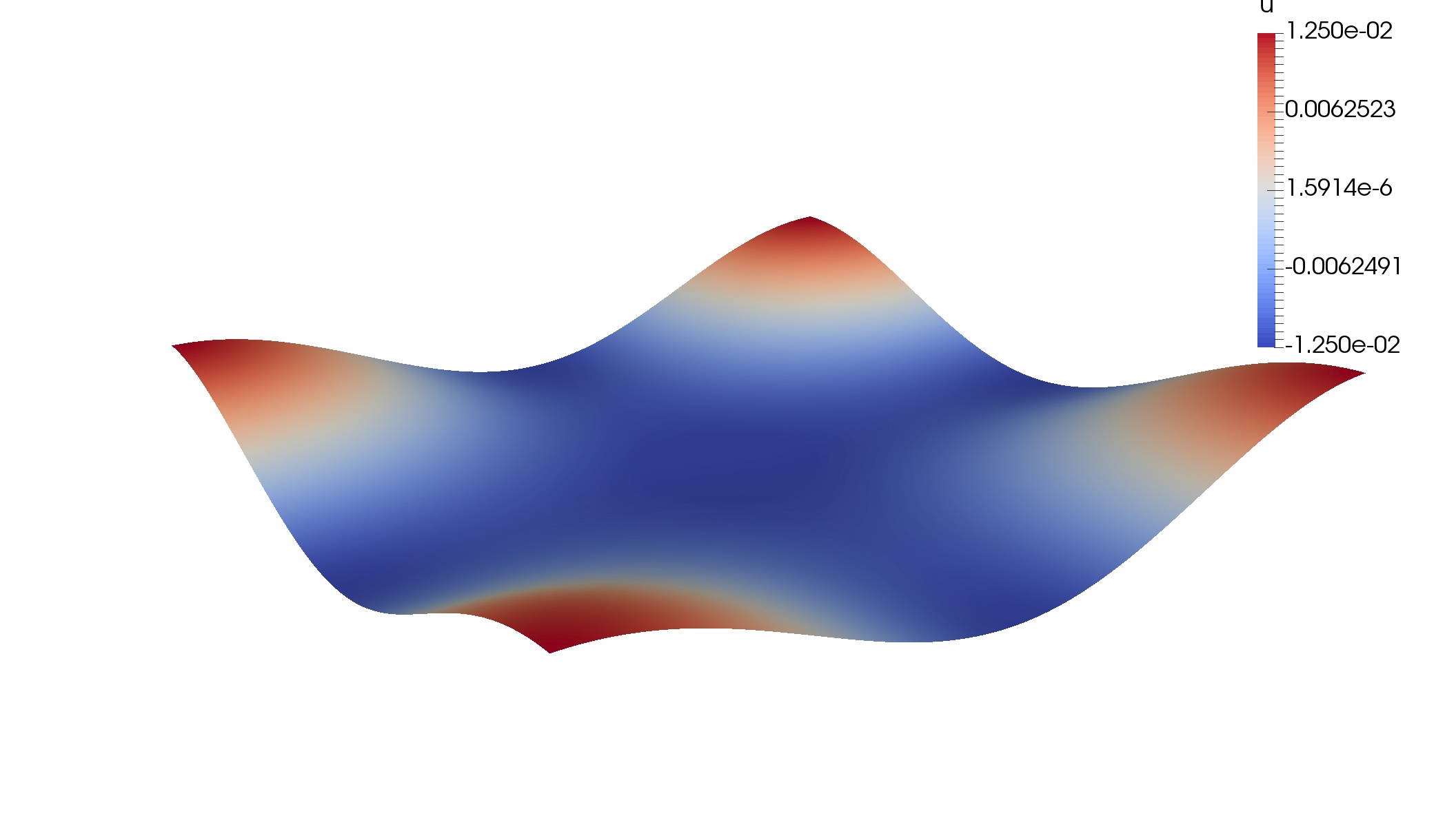}
    }  
    \hfill
    \subfigure[{\label{fig:a2}
        The approximation to $u$, the solution of the $142$-Bilaplacian.
    }]{
      \includegraphics[scale=\figscale,width=0.47\figwidth]{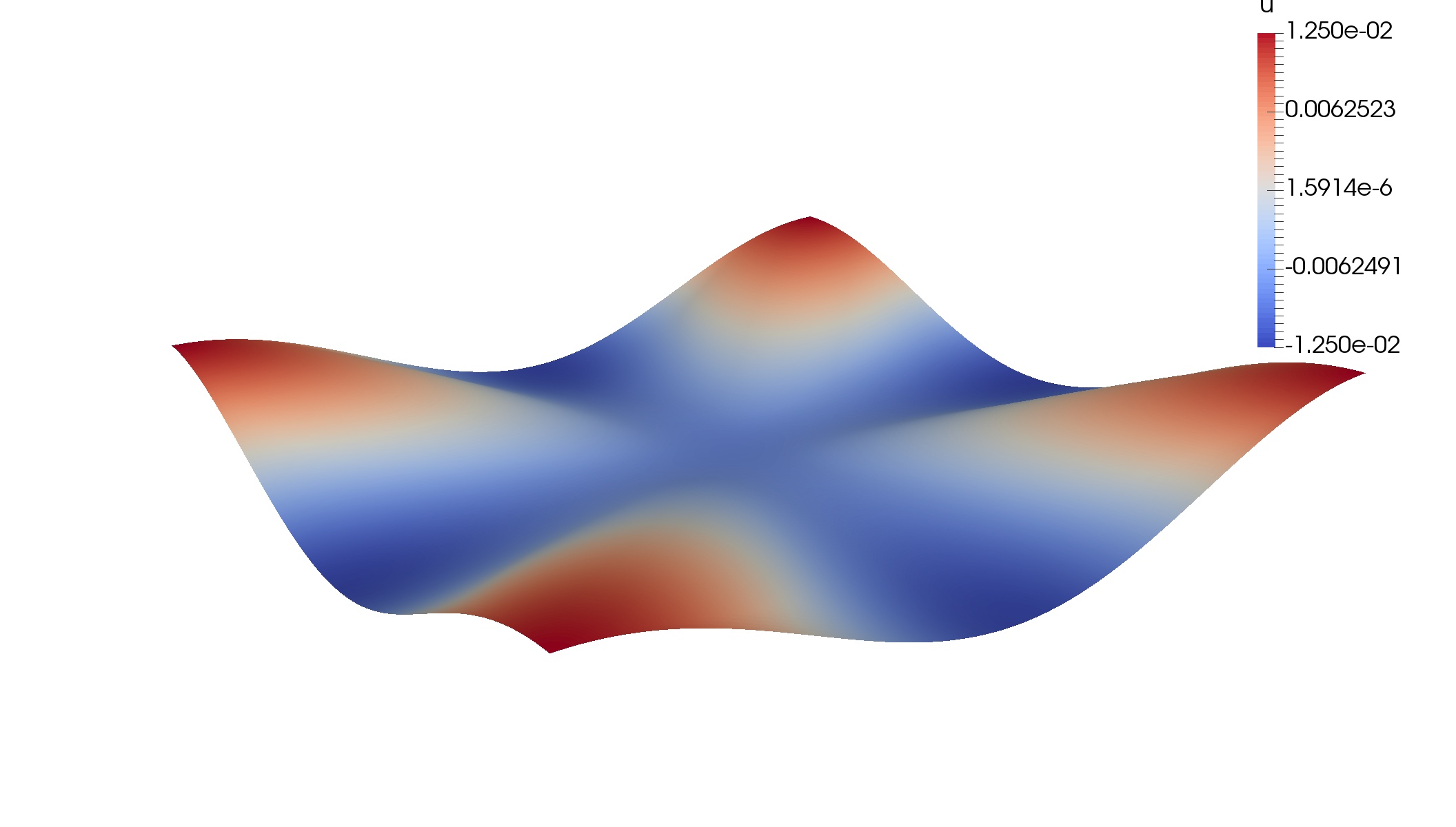}
    }  
  \end{center}
\end{figure}

}

\ms
\ms

 \ms

\ms
\noi \textbf{Acknowledgement.} N.K. would like to thank Craig Evans, Robert Jensen, Roger Moser, Juan Manfredi and Jan Kristensen for their inspiring mathematical discussions and especially their illuminating remarks on $\mD$-solutions and on $2$nd order $L^\infty$ variational problems.

\newpage

$ $

\newpage

\bibliographystyle{amsplain}

\end{document}